\documentclass[aos,seceqn,citesort,noautosecdot,dvips]{arximspdf}
\usepackage{graphics}
\doi{10.1214/09-AOS764}
\volume{38}
\issue{3}
\pubyear{2010}
\firstpage{1686}
\lastpage{1732}

\makeatletter

\newtheorem{theorem}{Theorem}[section]
\newtheorem{lemma}{Lemma}[section]

\makeatother

\begin{document}
\begin{frontmatter}

\title{Innovated higher criticism for detecting sparse signals in
correlated noise}
\runtitle{Innovated higher criticism}

\begin{aug}
\author[A]{\fnms{Peter} \snm{Hall}} and
\author[B]{\fnms{Jiashun} \snm{Jin}\corref{}\thanksref{t1}\ead[label=e1]{jiashun@stat.cmu.edu}}
\runauthor{P. Hall and J. Jin}
\affiliation{University of Melbourne and University of California at
Davis, and~Carnegie~Mellon~University}
\address[A]{Department of Mathematics and Statistics\\
University of Melbourne\\
Parkville, VIC, 3010\\
Australia\\
and\\
Department of Statistics\\
University of California\\
Davis, California 95616\\
USA} 
\address[B]{Department of Statistics\\
Carnegie Mellon University\\
Pittsburgh, Pennsylvania 15213\\
USA\\
\printead{e1}}
\end{aug}

\thankstext{t1}{Supported in part by NSF CAREER Award DMS-0908613.}

\received{\smonth{8} \syear{2009}}
\revised{\smonth{10} \syear{2009}}

%
\begin{abstract}
Higher criticism is a method for detecting signals that are both sparse
and weak. Although first proposed in cases where the noise variables
are independent, higher criticism also has reasonable performance in
settings where those variables are correlated. In this paper we show
that, by exploiting the nature of the correlation, performance can be
improved by using a modified approach which exploits the potential
advantages that correlation has to offer. Indeed, it turns out that
the case of independent noise is the most difficult of all, from a
statistical viewpoint, and that more accurate signal detection (for a
given level of signal sparsity and strength) can be obtained when
correlation is present. We characterize the advantages of correlation
by showing how to incorporate them into the definition of an optimal
detection boundary. The boundary has particularly attractive
properties when correlation decays at a polynomial rate or the
correlation matrix is Toeplitz.
\end{abstract}

%

\end{frontmatter}

\section{\texorpdfstring{Introduction.}{Introduction}} \label{sec:INT}

Donoho and Jin \cite{DJ04} developed Tukey's \cite{Tukey} proposal for
``higher criticism'' (HC), showing that a method based on the
statistical significance of a large number of statistically significant
test results could be used very effectively to detect the presence of
very sparsely distributed signals. They demonstrated that HC is capable
of optimally detecting the presence of signals that are so weak and so
sparse that the signal cannot be consistently estimated.
Applications include the problem of signal detection against cosmic
microwave background radiation (Cayon, Jin and Treaster \cite{Cayon2},
Cruz et al.
\cite{Cruz}, Jin \cite{Jin04,Jin06,Jin07}, Jin et al. \cite
{Starck}). Related work includes that of Cai, Jin and Low \cite{CJL},
Hall, Pittelkow and Ghosh
\cite{HPG} and Meinshausen and Rice \cite{Rice}.

The context of Donoho and Jin's \cite{DJ04} work was that where the
noise is white, although a small number of investigations have been
made of the case of correlated noise (Hall, Pittelkow and Ghosh \cite
{HPG}, Hall and
Jin \cite{HJ08}, Delaigle and Hall \cite{DH}). However, that research
has focused on the ability of standard HC, applied in the form that is
appropriate for independent data, to accommodate the nonindependent
case. In this paper we address the problem of how to modify HC by
developing \textit{innovated higher criticism} (iHC) and showing how to
optimize performance for correlated noise.

Curiously, it turns out that when using the iHC method tuned to give
optimal performance, the case of independence is the most difficult of
all, statistically speaking. To appreciate why this result is
reasonable, note that if the noise is correlated then it does not vary
so much from one location to a nearby location, and so is a little
easier to identify. In an extreme case, if the noise is perfectly
correlated at different locations then it is constant, and in this
instance it can be easily removed.

On the other hand, standard HC does not perform well in the case of
correlated noise, because it utilizes only the marginal information in
the data without much attention to the correlation structure. Innovated
HC is designed to exploit the advantages offered by correlation and
gives good performance across a wide range of settings.

The concept of the ``detection boundary'' was introduced by Donoho and
Jin \cite{DJ04} in the context of white noise. In this paper, we extend
it to the correlated case. In brief, the detection boundary describes
the relationship between signal sparsity and signal strength that
characterizes the boundary between cases where the signal can be
detected and cases where it cannot. In the setting of dependent data,
this watershed depends on the correlation structure of the noise as
well as on the sparsity and strength of the signal. When correlation
decays at a polynomial rate we are able to characterize the detection
boundary quite precisely. In particular, we show how to construct
concise lower/upper bounds to the detection boundary, based on the
diagonal components of the inverse of the correlation matrix, $\Sigma
_n$. A~special case is where $\Sigma_n$ is Toeplitz; there the upper
and the lower bounds to the detection boundary are asymptotically the
same. In the Toeplitz case, the iHC is optimal for signal detection but
standard HC is not.

There is a particularly extensive literature on multiple hypothesis
testing under conditions of dependence. It includes contributions to
the control of family-wise error rate and false discovery rate, and
work of Abramovich et al. \cite{ABDJ}, Benjamini and Hochberg \cite
{BenHoch}, Benjamini and Yekutieli \cite{BenYek}, Brown and Russel
\cite
{Brown}, Cai and Sun~\cite{CaiSun}, Clarke and Hall \cite{Clarke},
Cohen, Sackrowitz and Xu \cite{Cohen}, Donoho and Jin~\cite{DJ06},
Dunnett and
Tamhane \cite{Dunn}, Efron \cite{Efron}, Finner and Roters \cite{Finn},
Genovese and Wasserman \cite{Gen}, Jin and Cai \cite{JC}, Olejnik et al. \cite{Ole}, Rom
\cite
{Rom}, Sarkar and Chang \cite{Sar} and Wu \cite{Wu}. Work of Kuelbs and
Vidyashankar \cite{Kuelbs} is also related. Our contributions differ
from those of these authors in that we point to the advantages, rather
than the disadvantages, of dependence, and show how the advantages can
be exploited. In particular, as noted above, the problem of denoising
dependent data is actually simpler than in the case of independence. We
show how to exploit dependence and obtain improvements in performance
relative to what is possible in the context of independence and also
relative to the inferior performance that is obtained if a method that
is designed for the case of independence is applied inappropriately to
dependent data. In contrast, earlier work has tended to try to minimize
the problems caused by dependence rather than to capitalise on the
advantages that are available.

The paper is organized as follows. Section \ref{sec:review} introduces
the sparse signal model followed by a brief review of the uncorrelated
case. Section \ref{sec:LB} establishes lower bounds to the detection
boundary in correlated settings. Section \ref{sec:iHC} introduces
innovated HC and establishes an upper bound to the detection boundary.
Section \ref{sec:Toeplitz} applies the main results in Sections \ref
{sec:LB} and \ref{sec:iHC} to the case where the $\Sigma_n$'s are
Toeplitz. In this case, the lower bound coincides with the upper bound
and innovated HC is optimal for detection. Section \ref{sec:cluster}
discusses a case where the signals have a more complicated structure.
Section \ref{sec:strong} investigates a case of strong dependence.
Simulations are given in Section \ref{sec:Simul}, and discussion is
given in Section \ref{sec:Discu}. Section \ref{sec:proof}
and the \hyperref[sec:appen]{Appendix} give proofs of theorems and lemmas, respectively.

\section{\texorpdfstring{Sparse signal model and review of
HC.}{Sparse signal model and review of HC}} \label{sec:review}

\subsection{\texorpdfstring{Model.}{Model}}

Consider an $n$-dimensional Gaussian vector,
%
%
\begin{equation} \label{Model0}
X = \mu+ Z \qquad\mbox{where } Z \sim\mathrm{N}(0, \Sigma),
\end{equation}
with the mean vector $\mu$ unknown and the dimension $n$ large.
In most parts of the paper, we assume that $\Sigma= \Sigma_n$ is
known and has unit diagonal elements (the case where $\Sigma_n$ is
unknown is discussed in Section \ref{estimatingsigma} and Section \ref{sec:Discu}).
We are interested in testing whether no signal exists (i.e., $\mu= 0$)
or there is a sparse and faint signal.

Formulae (\ref{Modelpara1}) and (\ref{sigstrength}), below, introduce quantities $m$ and $A_n$
that represent signal sparsity and signal strength, respectively. In
particular, as $m$ increases the amount of sparsity decreases, and as
$A_n$ increases the strength of the signal increases. Of course, an
increase in either $m$ or $A_n$ leads to an increase in the ease with
which the signal can be detected and read. It would be possible to
connect $m$ and $A_n$ by a formula, and use that relationship to adjust
the signal, but we feel that the influence of the key elements of
sparsity and strength are most clearly presented by treating them
separately. In particular, we model the number of nonzero entries of
$\mu$ as
%
%
\begin{equation} \label{Modelpara1}
m = n^{1-\beta}\qquad \mbox{where $\beta\in(1/2, 1)$}.
\end{equation}
This is a very sparse case, for the proportion of signals is much
smaller than $1/\sqrt{n}$. We suppose that the signals appear at $m$ different
locations---$\ell_1 < \ell_2 < \cdots< \ell_m$---that are randomly
drawn from $\{1, 2, \ldots, n\}$ without replacement,
%
%
\begin{eqnarray} \label{Modelpara2}
P\{\ell_1 = n_1, \ell_2 = n_2, \ldots, \ell_m = n_m\} = \pmatrix
{n\cr
m}^{-1}\nonumber\\[-8pt]\\[-8pt]
\eqntext{\mbox{for all } 1 \leq n_1 < n_2 < \cdots< n_m \leq
n ,}
\end{eqnarray}
and that they have a common magnitude of
%
%
\begin{equation} \label{sigstrength}
A_n = \sqrt{2 r \log n} \qquad\mbox{where } r \in(0,1).
\end{equation}
These assumptions are made throughout the paper, in cases where $\Sigma
$ is relatively general as well as in cases (see Sections \ref
{sec:boundary} and \ref{sec:OHC}, below) where the noise variables are
assumed uncorrelated and so $\Sigma$ is the identity. Variations of
this model give similar results. For example, if we take the $j$th
nonzero signal to equal $W_j \sqrt{2\log n}$, where the $W_j$'s are
independent random variables with a common, nonnegative distribution
that has an upper endpoint $r^{1/2}$ satisfying $P(W\leq r^{1/2})=1$
and $P(W>r^{1/2}-\varepsilon)>0$ for all $\varepsilon>0$, then the results
are identical to their counterparts when signal strength is given
by (\ref{sigstrength}).

We are interested in testing which of the following two hypotheses is true:
%
%
\begin{equation} \label{ModelTest}
H_0\dvtx\mu= 0 \quad\mbox{vs.}\quad H_1^{(n)}\dvtx
\mu\mbox{ is a sparse vector as above}.
\end{equation}
This testing problem was found to be delicate even in the uncorrelated
case where $\Sigma_n = I_n$. See \cite{DJ04} (also \cite
{CJL,Ingster97,Ingster99,Jin04,Rice}) for details.

The case where $\Sigma_n$ is not the identity can arise when signals
are recorded at time points that are closely spaced in time or space.
See Section \ref{estimatingsigma} for discussion. An example of a
different type is that of \textit{global testing} in linear models. Here
we consider a model $Y \sim\mathrm{N}(M \mu, I_n)$, where the matrix
$M$ has
many rows and columns, and we are interested in testing whether $\mu=
0$. The setting is closely related to model (\ref{Model0}), since the
least squares estimator of $\mu$ is distributed as $\mathrm{N}(\mu, (M'
M)^{-1})$. The global testing problem is important in many
applications. One is that of testing whether a clinical outcome is
associated with the expression pattern of a pre-specified group of
genes (Goeman et al. \cite{Goeman1,Goeman2}) where $M$ is the
expression profile of the specified group of genes. Another is
expression quantitative Trait Loci (eQTL) analysis where $M$ is related
to the numbers of common alleles for different genetic markers and
individuals (Chen, Tong and Zhao \cite{Hongyu1}). In both examples,
$M$ is
either observable or can be estimated. Also, it is frequently seen that
only a small proportion of genes is associated with the clinical
outcome, and each gene contributes weakly to the clinical outcome. In
such a situation, the signals are both sparse and faint.


\subsection{\texorpdfstring{Detection boundary in the uncorrelated
case $(\Sigma_n =
I_n)$.}{Detection boundary in the uncorrelated case $(\Sigma_n =
I_n)$}} \label{sec:boundary}
\begin{figure}

\includegraphics{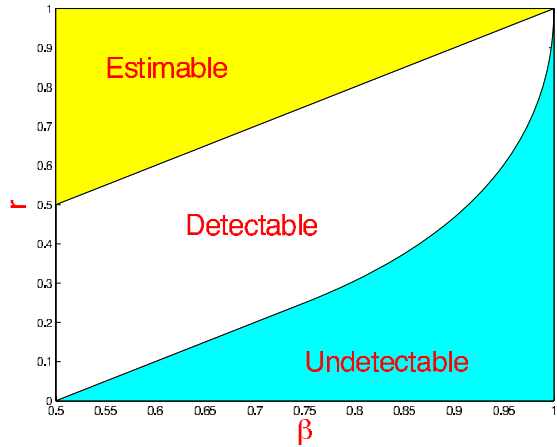}

\caption{Phase diagram for the detection problem in the uncorrelated
case. The detection boundary separates the $\beta$--$r$ plane into the
detectable region and the undetectable region. In the estimable region,
it is not only possible to reliably tell the existence of nonzero
coordinates, but is also possible to identify them individually. }
\label{fig:Detect}
\end{figure}

The testing problem is characterized by the curve $r = \rho^*(\beta)$
in the $\beta$--$r$ plane where
%
%
\begin{equation} \label{Definerho*1}
\rho^*(\beta)
= \cases{
\beta- 1/2, &\quad $1/2 < \beta\leq3/4$, \cr
\bigl(1 - \sqrt{1 - \beta}\bigr)^2, &\quad $3/4 < \beta< 1$,}
\end{equation}
and we call $r = \rho^*(\beta)$ the \textit{detection boundary}. The
detection boundary partitions the $\beta$--$r$ plane into two
sub-regions: the \textit{undetectable region} below the boundary and the
\textit{detectable region} above the boundary (see Figure \ref
{fig:Detect}). In the interior of the undetectable region, the signals
are so sparse and so faint that no test is able to successfully
separate the alternative hypothesis from the null hypothesis in (\ref
{ModelTest}): the sum of types I and II errors of any test tends to
$1$ as $n$ diverges to infinity. In the interior of the detectable
region, it is possible to have a test such that as $n$ diverges to
infinity, the type I error tends to zero and the power tends to $1$.
[In fact, Neyman--Pearson's Likelihood Ratio Test (LRT) is such a
test.]
See \cite{DJ04,Ingster97,Jin04}, for example.

%
%

The drawback of LRT is that it needs detailed information about the
unknown parameters $(\beta, r)$. In practice, we need a test that does
not need such information; this is where HC comes in.

\subsection{\texorpdfstring{Higher criticism and its optimal
adaptivity in the
uncorrelated case $(\Sigma_n = I_n)$.}{Higher criticism and its
optimal adaptivity in the
uncorrelated case $(\Sigma_n = I_n)$}} \label{sec:OHC}
A notion that goes back to Tukey \cite{Tukey}, higher criticism was
first proposed in \cite{DJ04} to tackle the aforementioned testing
problem in the uncorrelated case. To apply higher criticism, let $p_j =
P\{|\mathrm{N}(0,1)| \geq|X_j| \}$ be the $p$-value associated with
the $j$th
observation unit, and let $p_{(j)}$ be the $j$th $p$-value after
sorting in ascending order. The higher criticism statistic is defined as
%
%
\begin{equation} \label{DefineOHC}
\mathrm{HC}_n^* = \max_{j\dvtx1/n \leq p_{(j)} \leq1/2} \biggl\{
\frac
{\sqrt{n} (j/n - p_{(j)})}{\sqrt{p_{(j)} (1 - p_{(j)}) }} \biggr\}.
\end{equation}
There are also other versions of HC (see, e.g., \cite{DJ04,DJ08a,DJ08b}).
When $H_0$ is true, $\mathrm{HC}_n^*$ equals in distribution to the
maximum of the standardized uniform stochastic process \cite{DJ04}.
Therefore, by a well-known result for empirical processes \cite{Wellnerbook},
%
%
\begin{equation} \label{HCTest0}
\frac{\mathrm{HC}_n^*}{\sqrt{2 \log\log n}} \rightarrow1 \qquad\mbox{in
probability}.
\end{equation}
Consider the higher criticism test which rejects the null hypothesis when
%
%
\begin{equation} \label{HCTest}
\mathrm{HC}_n^* \geq(1 + a) \sqrt{2 \log\log n} \qquad\mbox{where
$a >
0$ is a constant}.
\end{equation}
It follows from (\ref{HCTest0}) that
the type I error tends to zero as $n$ diverges to infinity. For any
parameters $(\beta, r)$ that fall in the interior of the detectable
region, the type II error also tends to zero. This is the following theorem.
\begin{theorem} \label{thm:OHC}
Consider the higher criticism test that rejects $H_0$ when
$\mathrm{HC}_n^* \geq(1+a) \sqrt{2 \log\log n}$. For every alternative
$H_1^{(n)}$ where the associated parameters $(r, \beta)$ satisfy $r
> \rho^*(\beta)$, the HC test has asymptotically full power for
detection:
\[
P_{H_1^{(n)}} \{\mbox{Reject $H_0$}\} \rightarrow1 \qquad\mbox{as } n
\rightarrow\infty.
\]
\end{theorem}

That is, the higher criticism test adapts to unknown parameters $(\beta
,r)$ and yields asymptotically full power for detection throughout the
entire detectable region. We call this the \textit{optimal adaptivity} of
higher criticism \cite{DJ04}.

Theorem \ref{thm:OHC} is closely related to \cite{DJ04}, Theorem 1.2,
where a mixture model is used. The mixture model reduces approximately
to the current model if we randomly shuffle the coordinates of $X$.
However, despite its appealing technical convenience, it is not clear
how to
generalize the mixture model from the uncorrelated case to general
correlated settings. Theorem \ref{thm:OHC} is a special case of Theorem
\ref{thm:HC}.

We now turn to the correlated case. In this case, the exact ``detection
boundary'' may depend on $\Sigma_n$ in a complicated manner, but it is
possible to establish both a tight lower bound and a tight upper bound.
We discuss the lower bound first.

\section{\texorpdfstring{Lower bound to detectability.}{Lower bound
to detectability}} \label{sec:LB}

To establish the lower bound, a key element is the theory in comparison
of experiments (e.g., \cite{Strasser}) where a useful guideline is that
adding noise always makes the inference more difficult. Thus we can
alter the model by either adding or subtracting a certain amount of
noise so that the difficulty level (measured by the Hellinger distance,
or the $\chi^2$-distance, etc., between the null density and the
alternative density) of the original problem is sandwiched by those of
the two adjusted models. The correlation matrices in the latter have a
simpler form and hence are much easier to analyze. Another key element
is the recent development of matrix characterizations based on
polynomial off-diagonal decay where it shows that the inverse of a
matrix with this property shares the same rate of decay as the original matrix.

\subsection{\texorpdfstring{Comparison of experiments: Adding noise
makes inference
harder.}{Comparison of experiments: Adding noise makes inference harder}}
We begin by comparing two experiments that have the same mean, but
where the data from one experiment are more noisy than those from the
other. Intuitively, it is more difficult to make inference in the first
experiment than in the other. Specifically, consider the two Gaussian models
%
%
\begin{eqnarray} \label{CompareExperiment}
X &=& \mu+ Z,\qquad Z \sim\mathrm{N}(0, \Sigma)
\quad\mbox{and}\nonumber\\[-8pt]\\[-8pt]
X^* &=& \mu+ Z^*,\qquad Z^* \sim\mathrm{N}(0, \Sigma^*),\nonumber
\end{eqnarray}
where $\mu$ is an $n$-vector that is generated according to some
distribution $G = G_n$. The second model is more noisy than the first,
in the sense that $\Sigma^* \geq\Sigma$. Here, given two matrices, $A$
and $B$, we write $A \geq B$ if $A - B$ is positive semi-definite.

The second model in (\ref{CompareExperiment}) can be viewed as the
result of adding noise to the first. Indeed, defining $\Delta= \Sigma
^* - \Sigma$, taking $\xi$ to be N$(0, \Delta)$ (independently of $Z$),
and noting that $Z + \xi\sim\mathrm{N}(0, \Sigma+ \Delta)$, the
second model
is seen to be equivalent to $X + \xi= \mu+ (Z + \xi)$. Intuitively,
adding noise makes inference more difficult because it reduces the
distance between between $X$ and $Z$. To make this point concisely, let
$\operatorname{Hd}(X,Z;\mu,\Sigma)$ and $\operatorname
{Hd}(X^*,Z^*;\mu,\Sigma^*)$ denote
the Hellinger distance between (the distributions of) $X$ and $Z$, and
between $X^*$ and $Z^*$, respectively. Then we claim that the first of
these distances exceeds the second
%
%
\begin{equation} \label{monotone}\quad
\mbox{if } \Sigma^*\geq\Sigma\mbox{ in (\ref
{CompareExperiment})}\qquad
\mbox{then }
\operatorname{Hd}(X,Z;\mu,\Sigma)\geq\operatorname{Hd}(X^*,
Z^*;\mu,\Sigma^*).
\end{equation}
See Section \ref{sec:proof} for a proof. [The Hellinger distance
between distributions with densities $f$ and $g$ equals $\frac
{1}{2}\int
(f^{1/2}-g^{1/2})^2$.]

\subsection{\texorpdfstring{Matrices having polynomial off-diagonal
decay.}{Matrices having polynomial off-diagonal decay}}
\label{subsec:decayreview}
Next, we review results concerning matrices with polynomial
off-diagonal decay. The main message is that, under mild conditions, if
a matrix has polynomial off-diagonal decay, then its inverse as well as
its Cholesky factorization (which is unique if we require the diagonal
entries to be positive) also have polynomial off-diagonal decay, and
with the same rate. This beautiful result was recently obtained by
Jaffard \cite{Jaffard} (see also \cite{Grochenig1,Sun1}).

In detail, writing $\Theta_n$ for the set of $n \times n$ correlation
matrices, we introduce, for $\lambda>1$,
%
%
\begin{equation} \label{Assumption12}\hspace*{28pt}
\Theta_n^*(\lambda,c_0, M) = \{ \Sigma_n \in\Theta_n\dvtx
|\Sigma
_n(j,k)| \leq M(1 + |j - k|)^{-\lambda}, \|\Sigma_n\| \geq
c_0
\}.
\end{equation}
This is the set of matrices which have a given rate of polynomial
off-diagonal decay and where the operator norm is uniformly bounded
from below.
Consider a sequence of matrices $\{\Sigma_n\}_{n = 1}^{\infty}$ such
that $\Sigma_n \in\Theta_n^*(\lambda, c_0, M)$ for each $n$. It turns
out that the inverses (as well as the Cholesky factorizations) of
such sequences enjoy polynomial off-diagonal decay with the same rate
as that of the matrices themselves. See the \hyperref[sec:appen]{Appendix} for
the proof.

We are now ready for the lower bound.

\subsection{\texorpdfstring{Lower bound to detectability.}{Lower
bound to detectability}} \label{subsec:LB}
Consider a sequence of matrices $\{\Sigma_n\}_{n = 1}^{\infty}$ such
that $\Sigma_n \in\Theta_n^*(\lambda, c_0, M)$ for each $n$. Suppose
the extreme diagonal entries of $\Sigma_n^{-1}$ have an upper limit
$\bar{\gamma}_0$ in the range $0 < \bar{\gamma}_0 < \infty$; that is,
%
%
\begin{equation} \label{Assumption3}
\mathop{\overline{\lim}}_{n \rightarrow\infty} \Bigl( \max_{\sqrt{n}
\leq k \leq n -
\sqrt
{n} } \Sigma_n^{-1}(k,k) \Bigr) = \bar{\gamma}_0.
\end{equation}
Recall that the detection boundary in the uncorrelated case is defined
by $r<\rho^*(\beta)$. The following theorem asserts that, in the
presence of correlation, if we change the definition to $r<\bar{\gamma
}_0^{-1}\cdot\rho^*(\beta)$, then we obtain at least a lower bound to
the detection boundary.
\begin{theorem} \label{thm:LB}
Fix $\beta\in(1/2, 1)$, $r \in(0,1)$, $\lambda> 1$, $c_0 > 0$, and
$M > 0$. Consider a sequence of correlation matrices $\Sigma_n \in
\Theta_n^*(\lambda, c_0,M)$ that satisfy (\ref{Assumption3}). If
$r <
\bar{\gamma}_0^{-1} \rho^*(\beta)$, then the null hypothesis and
alternative hypothesis in (\ref{ModelTest}) merge asymptotically, and
the sum of types \textup{I} and \textup{II} errors of any test converges to $1$ as
$n$ diverges to infinity.
\end{theorem}

We now turn to the upper bound. The key is to adapt the higher
criticism to correlated noise and form a new statistic---\textit
{innovated higher criticism}.
\section{\texorpdfstring{Innovated higher criticism, upper bound to
detectability.}{Innovated higher criticism, upper bound to detectability}}
\label{sec:iHC}

Originally designed for the independent case, standard HC is
not really appropriate for dependent data for the following reasons.
First, HC only summarizes the
information that resides in the marginal effects of each coordinate and
neglects the correlation structure of the data.
Second, HC remains the same if we randomly shuffle different
coordinates of $X$. Such shuffling does not have an effect if $\Sigma_n
= I_n$, but does otherwise. In this section we build the correlation
into the standard higher criticism and form a new statistic---innovated
higher criticism (iHC). We then use
iHC to establish an upper bound to detectability. The iHC is intimately
connected to the well-known notion of innovation in time series \cite{Brock}
[see (\ref{ModelTransform1}) below],
hence the name innovated higher criticism.

Below, we begin by discussing the role of correlation in the detection problem.

\subsection{\texorpdfstring{Correlation among different coordinates:
Curse or
blessing\textup{?}}{Correlation among different coordinates: Curse or
blessing}}
\label{subsec:blessing}
Consider model (\ref{Model0}) in the two cases $\Sigma_n = I_n$ and
$\Sigma_n \neq I_n$. Which is the more difficult detection problem?

Here is one way to look at it. Since the mean vectors are the same in
the two cases, the problem where the noise vector contains more
``uncertainty'' is more difficult than the other. In information theory,
the \textit{total amount of uncertainty}
is measured by the \textit{differential entropy}, which in the Gaussian
case is proportional to the determinant of the correlation matrix \cite
{Cover}. As the determinant of a correlation matrix is largest when and
only when it is the identity matrix, the uncorrelated case contains the
largest amount of ``uncertainty'' and therefore gives the most difficult
detection problem. In a sense, the correlation is a ``blessing'' rather
than a ``curse'' as one might have expected.

Here is another way to look at it.
For any positive definite matrix $\Sigma_n$, denote the inverse of its
Cholesky factorization by $U_n$, a function of $\Sigma_n$ (so that $U_n
\Sigma_n U_n' = I_n$). Model (\ref{Model0}) is equivalent to
%
%
\begin{equation}\label{ModelTransform1}
U_n X = U_n \mu+ U_n Z \qquad\mbox{where } U_n Z \sim\mathrm{N}(0, I_n).
\end{equation}
(In the literature of time series \cite{Brock}, $U_n X$ is intimately
connected to the notion of innovation.) Compared to the uncorrelated
case, that is,
\[
X = \mu+ Z \qquad\mbox{where } Z \sim\mathrm{N}(0, I_n).
\]
It turns out that the noise vectors have the same distribution, but the
signals in the former are stronger. In fact, let $\ell_1 < \ell_2 <
\cdots< \ell_m$ be the $m$ locations where $\mu$ is nonzero. Recalling
that $\mu_j = A_n$ if $j \in\{ \ell_1, \ell_2, \ldots, \ell_m \}$,
$\mu_j = 0$ otherwise, and that $U_n$ is a lower triangular matrix,
%
%
\begin{equation} \label{ExplainUmu1}
(U_n \mu)_{\ell_k} = A_n \sum_{j = 1}^k U_n(\ell_k, \ell_i) = A_n
U_n(\ell_k, \ell_k) + A_n \sum_{j = 1}^{k - 1} U_n(\ell_j, \ell_k).
\end{equation}
Two key observations are as follows. First, since $\Sigma_n$ has unit
diagonal entries, every diagonal entry of $U_n$ is greater than or
equal to $1$, especially
%
%
\begin{equation} \label{ExplainUmu2}
U_n(\ell_k, \ell_k) \geq1.
\end{equation}
Second, recall that $m \ll n$, and $\{\ell_1, \ell_2, \ldots, \ell
_m \}
$ are randomly generated from $\{1, 2, \ldots, n\}$, so different
$\ell
_j$ are far apart from each other. Therefore, under mild decay
conditions on $U_n$,
%
%
\begin{equation} \label{ExplainUmu3}
U_n(\ell_j, \ell_k) \approx0,\qquad j = 1, 2, \ldots, k - 1.
\end{equation}
Inserting (\ref{ExplainUmu2}) and (\ref{ExplainUmu3}) into (\ref
{ExplainUmu1}),
we expect that $(U_n \mu)_{\ell_k} \gtrsim A_n$ for $k = 1, 2, \ldots,
m$. Therefore, ``on average,'' $U_n \mu$ has at least $m$ entries each
of which is at least as large as $A_n$.
This says that, first, the correlated case is easier for detection than
the uncorrelated case. Second, applying standard HC to $U_n X$ yields
a larger power than applying it to $X$ directly.

Next we make the argument more precise. Fix a positive sequence $\{
\delta_n\dvtx n \ge1 \} $ that tends to zero as $n$ diverges to infinity,
and a sequence of integers $\{b_n\dvtx n \geq1\}$ that satisfy $1 \leq b_n
\leq n$. Recall that $U_n$ is the function of $\Sigma_n$ defined by
\mbox{$U_n\Sigma_n U_n'=I_n$}, and let
\begin{eqnarray*}
&&\tilde{\Theta}_n^*(\delta_n, b_n) = \Biggl\{\Sigma_n \in\Theta_n,
\sum_{j = 1}^{k - b_n} |U_n(k,j)| \leq\delta_n,\\
&&\hspace*{69.1pt} \mbox{for all $k$
satisfying $ b_n + 1 \leq k \leq n$} \Biggr\}.
\end{eqnarray*}
Introducing $\tilde{\Theta}_n^*$ seems a digression from our original
plan of focusing on $\Theta_n^*$ (the set of matrices with polynomial
off-diagonal decay), but it is interesting in its own right. In fact,
compared to $\Theta_n^*$, $\tilde{\Theta}_n^*$ is much broader as it
does not impose much of a condition on $\Sigma_n(j,k)$ for $|j - k|
\leq b_n$. This helps to illustrate how broadly the aforementioned
phenomenon holds.
The following theorem is proved in Section~\ref{sec:proof}.
\begin{theorem} \label{thm:HCAdd}
Fix $\beta\in(1/2, 1)$ and $r \in(\rho^*(\beta), 1)$. Let $b_n =
n^{\beta}/3$, and let $\delta_n$ be a positive sequence that tends to
zero as $n$ diverges to infinity.
Suppose we apply standard higher criticism to $U_n X$ and we reject
$H_0$ if and only if the resulting score exceeds $(1+a) \sqrt{2 \log
\log n}$ where $a>0$. Then, uniformly in all sequences of $\Sigma_n$
satisfying $\Sigma_n \in\tilde{\Theta}_n^*( \delta_n, b_n)$,
\[
P_{H_0} \{\mbox{Reject $H_0$} \} + P_{H_1^{(n)}} \{\mbox{Accept
$H_0$}\} \rightarrow0,\qquad n \rightarrow\infty.
\]
\end{theorem}

Generally, directly applying standard HC to $X$ does not yield the same
result (e.g., \cite{HJ08}).

\subsection{\texorpdfstring{Innovated higher criticism: Higher
criticism based on
innovations.}{Innovated higher criticism: Higher criticism based on
innovations}} \label{subsec:iHC}
We have learned that applying standard HC to $U_n X$ yields better
results than applying it to $X$ directly. Is this the best we can do?
No, there is
still space for improvement. In fact, HC applied to $U_n X$ is a
special case of innovated higher criticism to be elaborated in this
section. Innovated higher criticism is even more powerful in detection.

\begin{figure}[b]

\includegraphics{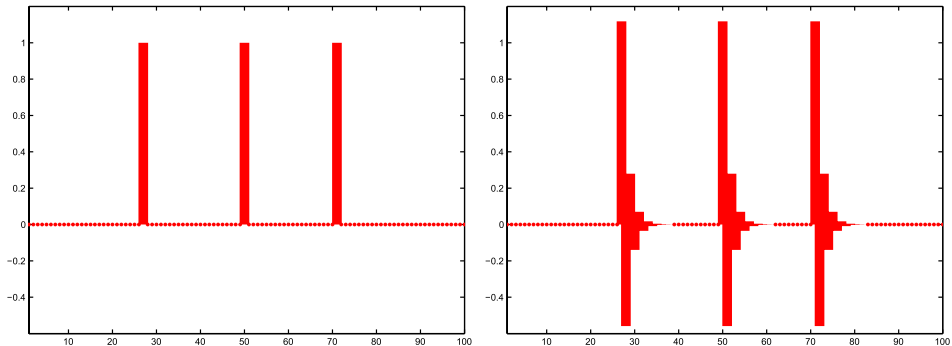}

\caption{Comparison of $\mu$ (left) and $U_n \mu$ (right). Here $n =
100$ and $\Sigma_n$ is a symmetric tri-diagonal matrix with $1$ on the
main diagonal, $0.4$ on two sub-diagonals and zero elsewhere. Also,
$\mu
$ is $1$ at coordinates $27$, $50$, and $71$ and $0$ elsewhere. In
comparison, the nonzero entries of $U_n \mu$ appear in three visible
clusters, each of which corresponds to a nonzero coordinate of $\mu$.}
\label{fig:Spikes}
\end{figure}

To begin, we revisit the vector $U_n \mu$ via an example.
Fix $n = 100$; let $\Sigma_n$ be a symmetric tri-diagonal matrix with
$1$ on the main diagonal, $0.4$ on two sub-diagonals and zero
elsewhere; and let $\mu$ be the vector with $1$ at coordinates $27$,
$50$, $71$ and zero elsewhere. Figure \ref{fig:Spikes} compares $\mu$
and $U_n \mu$. Especially, the nonzero coordinates of $U_n \mu$ appear
in three visible clusters, each of which corresponds to a different
nonzero entry of $\mu$. Also, at coordinates $27$, $50$, $71$, $U_n
\mu
$ approximately equals to $1.2$, but $\mu$ equals $1$. To interpret the
figure caption, recall that $U_n$ is the function of $\Sigma_n$ defined
by $U_n\Sigma_n U_n'=I_n$.

%
%

Now we can either simply apply standard HC to $U_n X$ as before, or we
can first linearly transform each cluster of signals to a singleton and
then apply the standard HC.
Note that in the second approach, we may have fewer signals, but each
of them is much stronger than those in $U_n X$.
Since the HC test is more sensitive to signal strength than to the
number of signals, we expect that the second approach yields greater
power for detection than the first.

In light of this we propose the following approach. Write $U_n =
(u_{kj})_{\{1 \leq k, j \leq n\}}$. We pick a bandwidth $1 \leq b_n
\leq n$, and construct a matrix
$\tilde{U}_n(b_n) = U_n(\Sigma_n, b_n)$ by banding $U_n$ \cite{Bickel}
%
%
\begin{equation} \label{DefineUbar}
\tilde{U}(b_n) \equiv(\tilde{u}_{kj} )_{1 \leq j, k \leq
n},\qquad
\tilde{u}_{kj} = \cases{
u_{kj}, &\quad $k - b_n + 1 \leq j \leq k$, \cr
0, &\quad otherwise.}
\end{equation}
We then normalize each column of $\tilde{U}_n(b_n)$ by its own $\ell
^2$-norm, and call the resulting matrix $\bar{U}_n(b_n)$. Next, defining
%
%
\begin{equation} \label{DefineVn}
V_n(b_n) =V_n(b_n; \Sigma_n) = \bar{U}_n' (b_n; \Sigma_n) \cdot U_n,
\end{equation}
we transform model (\ref{Model0}) into
%
%
\begin{equation} \label{ModelTransformx}
X \longmapsto V_n(b_n) X = V_n(b_n) \mu+ V_n(b_n) Z.
\end{equation}
Finally, we apply standard higher criticism to $V_n(b_n) X$, and call
the resulting statistic \textit{innovated higher criticism},
%
%
\begin{eqnarray}\label{DefineiHC}
\mathrm{iHC}_n^*(b_n) &=& \mathrm{iHC}_n^*(b_n;
\Sigma_n)\nonumber\\[-8pt]\\[-8pt]
&=& \frac
{1}{\sqrt
{2b_n - 1}} \sup_{j\dvtx1/n \leq p_{(j)} \leq1/2} \biggl\{ \sqrt{n}
\cdot\frac{j/n - p_{(j)}}{\sqrt{ p_{(j)} (1 - p_{(j)})} } \biggr\}.\nonumber
\end{eqnarray}
Note that standard HC applied to $U_nX$ is a special case of
$\mathrm{iHC}_n^*$ with $b_n = 1$.

We briefly comment on the selection of the bandwidth parameter $b_n$. First,
for each $k \in\{\ell_1, \ell_2, \ldots, \ell_m \}$, direct
calculations show that $(V_n(b_n) \mu)_k \approx A_n \cdot\sqrt{\sum
_{j = 1}^{b_n} u_{k, k - j + 1}^2} \geq A_n$. Second, $V_n(b_n) Z
\sim\mathrm{N}(0, \bar{U}_n'(b_n) \bar{U}_n(b_n))$, where
$\bar{U}_n'(b_n)\times
\bar
{U}_n(b_n)$ is a banded correlation matrix with bandwidth $2b_n-1$.
Therefore, choosing $b_n$ involves a trade-off: a larger $b_n$ usually
means stronger signals but also means stronger correlation among the
noise. While it is hard to give a general rule for selecting the best
$b_n$, we must mention that
in many cases, the choice of $b_n$ is not very critical. For example,
when $\Sigma_n$ has polynomial off-diagonal decay, a~logarithmically
large $b_n$ is usually appropriate.

\subsection{\texorpdfstring{Upper bound to detectability.}{Upper
bound to detectability}}
We now establish an upper bound to detectability. Suppose the diagonal
entries of
$\Sigma_n^{-1}$ have a lower limit as follows:
%
%
\begin{equation} \label{Assumption4}
\mathop{\underline{\lim}}_{n \rightarrow\infty} \Bigl( \min_{\sqrt{n}
\leq k \leq n -
\sqrt{n}}
\Sigma_n^{-1}(k,k) \Bigr) = \underline{\gamma_0}.
\end{equation}
Recall that the nonzero coordinates of $\mu$ are modeled as $A_n =
\sqrt
{2 r \log n}$. If we let $b_n = \log n$ then it can be proved that the
vector $V_n(b_n) \cdot X$ has at least $m$ nonzero coordinates, each of
which is as large as $\sqrt{\underline{\gamma_0}} A_n = \sqrt{2
\underline{\gamma_0} \cdot r \cdot\log n}$. (See Lemma~\ref
{lemma:Transformmean}.) Note that a larger $b_n$ cannot improve the
signal strength significantly, but may yield a much stronger
correlation in $V_n(b_n) Z$. Therefore, a smaller bandwidth is
preferred. The choice $b_n = \log n $ is mainly for convenience, and
can be modified.

We now turn to the behavior of $\mathrm{iHC}_n^*(b_n)$ under the null
hypothesis. In the independent case, $\mathrm{iHC}_n^*$ reduces to
$\mathrm{HC}_n^*$ and is approximately equal to $\sqrt{2 \log\log n}$. In the
current situation, $\mathrm{iHC}_n^*$ is comparably larger due to the
correlation. However, since the selected bandwidth is relatively small,
$\mathrm{iHC}_n^*$ remains logarithmically large. See Lemma \ref
{lemma:HCNull} for details.
The following theorem elaborates on the upper bound, and is proved in
Section \ref{sec:proof}.
\begin{theorem} \label{thm:HC}
Fix $c_0 > 0$, $\lambda> 1$, and $M > 0$, and set $b_n = \log n$.
Suppose $ \underline{\gamma_0} \cdot r > \rho^*(\beta)$. If we reject
$H_0$ when $\mathrm{iHC}_n^*(b_n; \Sigma_n) \geq(\log n)^{2}$, then,
uniformly in all $\Sigma_n \in\Theta_n^*(\lambda, c_0, M)$,
\[
P_{H_0} \{\mbox{Reject $H_0$} \} + P_{H_1^{(n)}} \{\mbox{Accept
$H_0$}\} \rightarrow0\qquad \mbox{as } n \rightarrow\infty.
\]
\end{theorem}

The cut-off value $(\log n)^2$ can be replaced by other logarithmically
large terms that tend to infinity faster than $(\log n)^{3/2}$. For
finite $n$, this cut-off value may be conservative. In Section \ref
{sec:Simul} [i.e., experiment (a)], we suggest an alternative where we
select the cut-off value by simulation.

In summary, a lower bound and an upper bound are established as $r =
\bar{\gamma}_0^{-1} \rho^*(\beta)$ and $r = \underline{\gamma_0}^{-1}
\rho^*(\beta)$, respectively, under reasonably weak off-diagonal decay
conditions. When
$\bar{\gamma}_0 = \underline{\gamma_0}$,
the gap between the two bounds disappears, and iHC is optimal for detection.
Below in Sections \ref{sec:Toeplitz}--\ref{sec:strong}, we investigate several Toeplitz cases, ranging
from weak dependence to strong dependence; for these cases, iHC is
optimal in detection.

\subsection{\texorpdfstring{Effect of estimating $\Sigma_n$ and related
issues.}{Effect of estimating $\Sigma_n$ and related issues}}
\label{estimatingsigma}
So far, we have assumed that the covariance matrix $\Sigma_n$ is known.
When $\Sigma_n$ is unknown, we could still use iHC if $\Sigma_n$ could
be estimated. We now briefly comment on
the effect of estimating $\Sigma_n$.

In practical problems where iHC methodology would be used, noise could
reasonably be represented as a time series, and its characteristics
estimated from data. In particular, the time series might be an
autoregression, and data over a longer period than that for which the
current dataset was recorded could be used to deduce properties of the
noise. Examples include detection of xenon byproducts as evidence of a
nuclear explosion, early detection of bioweapons and detection of
covert communications.

If data are gathered over a time period of length $p$, if the signal is
present at no more than $m=n^{1-\beta}$ points where $\beta\in(1/2,1)$
and if the maximum size of the signal is no greater than a constant
multiple of $(\log n)^{1/2}$, then it is typically possible to estimate
the components of $U_n$ at rate $(p^{-1} \log p)^{1/2}$ uniformly in
all components. From this property it can be proved that the difference
between $U_n\mu$ and its empirical form equals $O_p[\{n^{2-\beta}
p^{-1} (\log n) (\log p)\}^{1/2}]$, uniformly in all components.
Similarly, if the noise process is conventional (e.g., an
autoregression) then the distance between $U_nZ$ and its empirical form
can be shown to equal $O_p(n^{1+\varepsilon}/p)$ for all $\varepsilon>0$.
Therefore the effects of variance estimation will be asymptotically
negligible if, for some $\varepsilon>1-\beta$, $n^{1+\varepsilon
}/p\to0$
converges to zero as $n\to\infty$.

To appreciate the extent to which this condition is restrictive,
consider the case where the signals are particularly sparse, that
is, $\beta$ is close to 1; say, $\beta=1-\eta_1$ where $\eta_1>0$ is
small. Then the condition holds if $p$ is at least as large as
$n^{1+\eta_2}$ for some $\eta_2>\eta_1$. That is, the amount of time
for which data have to be acquired in order to estimate $\Sigma_n$ with
sufficient accuracy need only be a factor $n^\varepsilon$ greater than~$n$,
for $\varepsilon>0$ relatively small. As the prevalence of the signal
increased, the size of $\varepsilon$ would have to too.

Application of our methods to other problems, such as those involving
genomic data, can be inhibited by the difficulty of estimating $\Sigma
_n$ without information from outside the dataset. However, while there
is sometimes evidence of strong dependence in genomic data, from other
viewpoints the overall level of correlation is often quite low. For
example, Messer and Arndt \cite{Messer} argue that correlation decays
from about 0.08, at a separation of approximately two base pairs, to
about $0.01$ for a separation of ten base pairs. Work of Mansilla et
al. \cite{Mansilla} corroborates these figures. Results such as these,
together with the upper tail independence property which is generally
available for light-tailed distributions, suggest that for genomic data
it is possible to work effectively under the assumption that expression
levels are statistically independent, even when they are not. Details
are given by Delaigle and Hall \cite{DH}, who use the fact that in the
case of genomic data the variables are typically $t$-statistics.

More generally, cases where the signals are distributed nonrandomly
can be compared readily with the case of independent,
randomly-distributed signals, noted just below (\ref{Modelpara1}), as follows. Let
us take as our benchmark the classical problem $\mathcal{P}(n_0,m_0)$
where there are $n_0$ independent noise variables, and $m_0$ signals
are distributed randomly among the $n_0$ locations. We shall compare it
with the more general problem where the noise variables are
$d$-dependent with the integer $d$ depending on $n$. To quantify the
effects of nonrandomness we assume that $m=n^{1-\beta}$ signals are
distributed among $m/K$ clumps of length $K=K(n)$, and that the points
in clumps that are furthest to the right are distributed sequentially
among the integers $K,K+1,\ldots,n-1,n$, with each placement being
conditional on the clump not overlapping any pre-existing clumps. We
make no other assumption about the dependence structure of the process
for placing the clumps, only that it be independent of the noise
variables; and we assume that $d\leq K$. If $K=O(n^\eta)$ for all
$\eta
>0$ then, for each $\eta$, the difficulty of the signal detection
problem is bounded above by that of $\mathcal{P}(n^{1-\eta},mn^{-\eta
})$, and below by that of $\mathcal{P}(n,m)$. Since $\eta$ here is
arbitrary then it can be deduced that the effect of clustering has
asymptotically negligible effect. On the other hand, if $K=n^\eta
\ell
_n$ for a fixed $\eta>0$ and a quantity $\ell_n$ that satisfies $\ell
_n=O(n^\varepsilon)$ and $n^\varepsilon=O(\ell_n)$ for all
$\varepsilon>0$, then
the problem can be asymptotically as difficult as $\mathcal
{P}(n^{1-\eta
},mn^{-\eta})$ for the given value of $\eta$.

\section{\texorpdfstring{Application in the Toeplitz
case.}{Application in the Toeplitz case}} \label{sec:Toeplitz}

In this section, we discuss the case where $\Sigma_n$ is a (truncated)
Toeplitz matrix that is generated by a spectral density $f$ defined
over $(-\pi, \pi)$. In detail,
let $a_k= (2\pi)^{-1} \int_{|\theta| < \pi}f(\theta) e^{ - i k
\theta
} \,d\theta$ be the $k$th Fourier coefficient of $f$. The $n$th
truncated Toeplitz matrix generated by $f$ is the matrix $\Sigma_n(f)$
of which the $(j,k)$th element is $a_{j-k}$, for $1 \leq j, k \leq n$.

We assume that $f$ is symmetric and positive, that is,
%
%
\begin{equation} \label{Assumption3b}
c_0(f) \equiv\mathop{\operatorname{essinf}}_{-\pi\leq\theta\leq\pi} f(\theta)
> 0.
\end{equation}
First, note that $f$ is a density, so $a_0 = 1$ and $\Sigma_n(f)$ has
unit diagonal entries. Second, from the symmetry of $f$, it can be seen
that $\Sigma_n(f)$ is a real-valued symmetric matrix. Last, it is well
known \cite{Bottcher} that the smallest eigenvalue of $\Sigma_n(f)$ is
no smaller than $c_0(f)$, so $\Sigma_n(f)$ is positive definite.
Putting all these together, $\Sigma_n(f)$ is seen to be a correlation matrix.

Toeplitz matrices enjoy convenient asymptotic properties. In detail,
let $\lambda>1$ and suppose that additionally $f$ has at least
$\lambda
$ bounded derivatives [meaning, if $\lambda$ is a positive integer,
that $|f^{(j)}|$ is bounded for $0\leq j\leq\lambda$, and, if
$\lambda$
is not\vspace*{1pt} an integer, that $|f^{(j)}|$ is bounded for $0\leq j<\lambda$
and $|f^{(\lambda')}(\theta_1)-f^{(\lambda')}(\theta_2)|/|\theta
_1-\theta_2|^{\lambda-\lambda'}$ is bounded, where $\lambda'$ denotes
the largest integer less than $\lambda$]. Then by elementary Fourier
analysis, there is a constant $M_0 = M_0(f) > 0$ such that
%
%
\begin{equation} \label{Assumption3a}
|a_k| \leq M_0(f) (1 + k)^{- \lambda} \qquad\mbox{for } k = 0,
1,2, \ldots.
\end{equation}
Comparing (\ref{Assumption3b}) and (\ref{Assumption3a}) with the
definition of $\Theta_n^*$, we conclude that
%
%
\begin{equation} \label{ToeplitzMain1}
\Sigma_n \in\Theta_n^*(\lambda, c_0(f), M_0(f)).
\end{equation}

In addition, it is known that the inverse of $\Sigma_n(f)$ is typically
asymptotically equivalent to the Toeplitz matrix generated by $1/f$.
The diagonal entries of $\Sigma_n(1/f)$ are the well-known \textit{Wiener
interpolation rates} \cite{Wiener},
%
%
\begin{equation} \label{DefineCf}
C(f) = \frac{1}{2\pi} \int_{-\pi}^{\pi} \frac{1}{f(\theta)}\,
d\theta.
\end{equation}
From this property and a result of \cite{Bottcher}, Theorem 2.15, it
can be proved that
\[
{\max_{ \sqrt{n} \leq k \leq n - \sqrt{n}}} |\Sigma
_n^{-1}(f)(k,k) -
C(f) | = o(1).
\]
Comparing this with (\ref{Assumption3}) and (\ref{Assumption4}) we
deuce that
%
%
\begin{equation} \label{ToeplitzMain2}
\bar{\gamma}_0 = \underline{\gamma_0} = C(f).
\end{equation}
Combining (\ref{ToeplitzMain1}) and (\ref{ToeplitzMain2}), the
following theorem is a direct result of Theorems~\ref{thm:LB} and \ref
{thm:HCAdd} (the proof is omitted).
\begin{theorem} \label{thm:Toeplitz}
Fix $\lambda> 1$, and let $\Sigma_n(f)$ be the Toeplitz matrix
generated by a symmetric spectral density $f$ that satisfies
(\ref{Assumption3b}) and (\ref{Assumption3a}).
When $C(f) \cdot r < \rho^*(\beta)$, the null and alternative
hypotheses merge asymptotically, and the sum of types \textup{I} and
\textup{II}
errors of any test converges to $1$ as $n$ diverges to infinity.
When $C(f) \cdot r > \rho^*(\beta)$, suppose we apply iHC with
bandwidth $b_n = \log n$ and reject the null hypothesis when $\mathrm
{iHC}_n^*(b_n, \Sigma_n(f)) \geq(\log n)^2$. Then the type \textup{I} error of
iHC converges to zero, and its power converges to 1.
\end{theorem}

\begin{figure}

\includegraphics{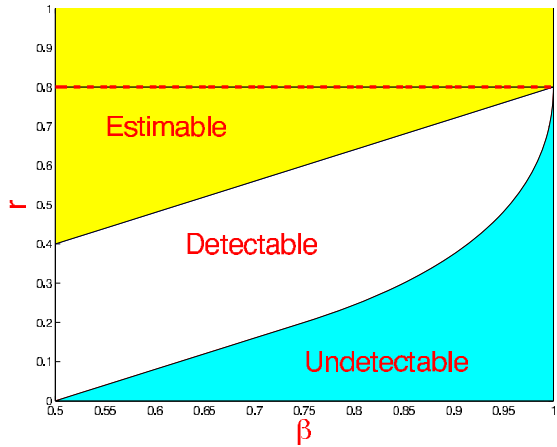}

\caption{Phase diagram in the case where $\Sigma_n$ is
a Toeplitz matrix generated by a spectral density~$f$. Similarly, in
Figure \protect\ref{fig:Detect}, the $\beta$--$r$ plane is
partitioned into three regions---undetectable, detectable,
estimable---each of which can be viewed as the corresponding region in
Figure \protect\ref{fig:Detect} squeezed vertically by a factor of $1/C(f)$.
In the rectangular region on the top, the largest signals in
$V_n(b_n) \cdot X$ [see (\protect\ref{DefineVn})] are large enough to stand
out by themselves. } \label{fig:DetectToep}
\end{figure}

The curve $r = C(f)^{-1} \rho^*(\beta)$ partitions the $\beta$--$r$
plane into the undetectable region and the detectable region, similarly
to the uncorrelated case. The regions of the current case can be viewed
as the corresponding regions in the uncorrelated squeezed vertically by
a factor of $1/C(f)$. See Figure \ref{fig:DetectToep}.
%
%
[Note that $C(f) \geq1$, with equality if and only if $f \equiv1$, which corresponds
to the uncorrelated case.]

\section{\texorpdfstring{Extension: When signals appear in
clusters.}{Extension: When signals appear in clusters}} \label{sec:cluster}

In the preceding sections [see, e.g., (\ref{Modelpara2}) in Section
\ref{sec:review}], the $m$ locations of signals were generated randomly
from $\{1, 2, \ldots, n\}$. Since $m \ll\sqrt{n}$, the signals appear
as singletons with overwhelming probabilities. In this section we
investigate an extension where the signals may appear in clusters.

We consider a setting where the signals appear in a total of $m$
clusters, whose locations are randomly generated from $\{1, 2, \ldots,
n\}$. Each cluster contains a total of $K$ consecutive signals, whose
strengths are $g_0 A_n$, $g_1 A_n, \ldots, g_{K-1} A_n$, from right
to left.
Here, $A_n = \sqrt{2 r \log n}$ as before, $K \geq1$ is a fixed
integer and $g_i$ are constants. Approximately, the signal vector can
be modeled as follows.

As before, let $\ell_1, \ell_2, \ldots, \ell_m$ be indices that are
randomly sampled from $\{1, 2, \ldots, n\}$. Let $\mu=(\mu_1,\ldots
,\mu
_n)^{\mathrm{T}}$, where $\mu_j= A_n$ if $j \in\{\ell_1, \ell_2,
\ldots
, \ell
_m\}$, and $\mu_j=0$ otherwise. Let $B = B_n$ denote the ``backward
shift'' matrix with 0 in every position except that it has 1 in position
$(j+1,j)$ for $1 \leq j \leq n-1$. Thus, $B\mu$ differs from $\mu$ in
that the components are shifted one position backward, with 0 added at
the bottom. We model the signal vector as
\[
\nu= g_0 \mu+ g_2 B \mu+ \cdots+ g_k B^{K-1} \mu
= \Biggl(\sum_{k = 0}^{K-1} g_k B^k \Biggr) \mu.
\]
Thus $\nu$ is comprised of $m$ clusters, each of which contains $K$
consecutive signals. Let $g$ be the function $g(\theta) = \sum_{0\leq
k\leq K-1} g_k e^{- i k \theta}$. We note that
$\sum_{0\leq k\leq K-1} g_k B^k$ is the lower triangular Toeplitz
matrix generated by $g$. With the same spectral density~$f$, we
consider an extension of that in Section \ref{sec:Toeplitz} by
considering the following model:
%
%
\begin{equation} \label{clustermodel}
X = \Sigma_n(g) \mu+ Z \qquad\mbox{where } Z \sim\mathrm{N}
(0, \Sigma_n(f)) ,
\end{equation}
with $f$ denoting the spectral density in Section \ref{sec:Toeplitz}.

We note that the model can be equivalently viewed as
\[
\tilde{X} = \mu+ \tilde{Z} \qquad\mbox{where } \tilde{Z}
\sim
\mathrm{N}(0, \tilde{\Sigma}_n) \quad\mbox{and}\quad \tilde{\Sigma}_n =
\Sigma_n^{-1}(g) \cdot\Sigma_n(f) \cdot\Sigma_n^{-1}(\bar{g}),
\]
with $\bar{g}$ denoting the complex conjugate of $g$. Asymptotically,
\[
\tilde{\Sigma}^{-1}_n \sim\Sigma_n(\bar{g}) \cdot\Sigma_n^{-1}(f)
\cdot\Sigma_n(g) \sim\Sigma_n(|g|^2/f),
\]
where the diagonal entries of $\Sigma_n(|g|^2/ f)$ are
\[
C(f, g) = \frac{1}{2 \pi} \int_{-\pi}^{\pi} \frac{|g(\theta
)|^2}{f(\theta)}\, d \theta.
\]
If $\bar{\gamma}_0$ and $\underline{\gamma_0}$ are as defined in
(\ref
{Assumption3}) and (\ref{Assumption4}), then
$\underline{\gamma_0} = \bar{\gamma}_0 = C(f, g)$,
and we expect the detection boundary to be
$r = C(f, g)^{-1} \cdot\rho^*(\beta)$.
This is affirmed by the following theorem which is proved in Section
\ref{sec:proof}.
\begin{theorem} \label{thm:cluster}
Fix $\lambda> 1$. Suppose $g_0 \neq0$ and let $f$ be a symmetric
spectral density that satisfies (\ref{Assumption3b}) and (\ref
{Assumption3a}). When $C(f, g) \cdot r < \rho^*(\beta)$, the null and
alternative hypotheses merge asymptotically, and the sum of types \textup{I} and
\textup{II} errors of any test converges to $1$ as $n$ diverges to infinity.
When $C(f, g) \cdot r > \rho^*(\beta)$, if we apply iHC to
$\Sigma_n^{-1} (g)X$ with bandwidth $b_n = \log n$ and reject the null
hypothesis when $\mathrm{iHC}_n^*(b_n, \Sigma_n^{-1}(g) \Sigma_n(f)
\Sigma
_n^{-1}(\bar{g})) \geq(\log n)^2$, then the type \textup{I} error converges to
zero, and the power converges to 1.
\end{theorem}

\section{\texorpdfstring{The case of strong dependence.}{The case of
strong dependence}} \label{sec:strong}

So far, we have only discussed weakly dependent cases. In this section,
we investigate the case of strong dependence.

Suppose we observe an $n$-variate Gaussian vector $X = \mu+ Z$,
where $\mu$ contains a total of $m$ signals, of equal strength to be
specified, whose locations are randomly drawn from $\{1,2, \ldots, n\}$
without replacement, and $Z \sim\mathrm{N}(0, \Sigma_n)$ where we
assume that
$\Sigma_n$ displays slowly decaying correlation,
%
%
\begin{equation} \label{hjmatrix}
\Sigma_n(j,k) = \max\{0, 1 - |j-k|^{\alpha} n^{-\alpha_0} \} ,\qquad
1 \leq j, k \leq n,
\end{equation}
with $\alpha> 0$ and $0 < \alpha_0 \leq\alpha$. The range of
dependence can be calibrated in terms of $k_0 = k_0(n; \alpha, \alpha
_0)$, denoting the largest integer by $k < n^{\alpha_0/\alpha}$.
Clearly, $k_0 \approx n^{\alpha_0/\alpha}$. Seemingly, the most
interesting range is $0 < \alpha_0\leq\alpha\leq1$.

Condition (\ref{hjmatrix}) is more restrictive than similar assumptions
in other places in this paper. There are at least two reasons. First,
the constants in the definition of the detection boundary turn out to
depend intimately on the value of $\alpha$ used in the definition of
$\Sigma_n$ at (\ref{hjmatrix}), and so we need to make an assumption
which is driven by that parameter. Secondly, a significantly more
general definition of $\Sigma_n$ would need to satisfy the positive
definiteness property which (as can be seen from Lemma~\ref{lemma:PD})
is somewhat delicate.


Model (\ref{hjmatrix}) has been studied in detail by Hall and Jin
\cite{HJ08} who showed that the detectability of standard HC is seriously
damaged by strong dependence. However, it remains open as to what is
the detection boundary, and how to adapt HC to overcome the strong
dependence and obtain optimal detection. This is what we address in the
current section.

The key idea is to decompose the correlation matrix as the product of
three matrices
each of which is relatively easy to handle. To begin with we introduce
a spectral density,
%
%
\begin{equation} \label{Definef0}
f_{\alpha}(\theta) = 1 - \sum_{k=1}^{\infty}
[(k+1)^{\alpha} + (k-1)^{\alpha} - 2 k^{\alpha} ] \cos
(k\theta) .
\end{equation}
[Note that the Fourier coefficients of $f_{\alpha}(\theta)$ satisfy the
decay condition in (\ref{Assumption3a}) with $\lambda= 2 - \alpha$.]
Next, let
\[
g_0(\theta)=1 - e^{-i\theta},\qquad a_n = a_n(\alpha_0) = n^{\alpha_0}/2.
\]
The Toeplitz matrix $\Sigma_n(g_0)$ is a lower triangular matrix with
$1$'s on the main diagonal, $-1$'s on the sub-diagonal and $0$'s elsewhere.
Additionally, let $D_n$ be the diagonal matrix where on the diagonal
the first entry is $1$ and the remaining entries are $\sqrt{a_n}$. Let
$\tilde{X} = D_n \cdot\Sigma_n(g_0) \cdot X$. Then model (\ref
{hjmatrix}) can be rewritten equivalently as
%
%
\begin{equation}\label{hj-model-1}
\tilde{X} = \tilde{\mu} + \tilde{Z} \qquad\mbox{where }
\tilde
{\mu} = D_n \cdot\Sigma_n(g_0) \cdot\mu\mbox{ and }
\tilde{Z} \sim\mathrm{N}(0, \tilde{\Sigma}_n)
\end{equation}
with $\tilde{\Sigma}_n = D_n \cdot\Sigma_n(g_0) \cdot\Sigma_n
\cdot\Sigma_n(\bar{g}_0) \cdot D_n$.
The key is that $\tilde{\Sigma}_n$ is asymptotically equivalent to the
Toeplitz matrix generated by $f_{\alpha}$. In detail, introduce
\[
\bar{\Sigma} = \pmatrix{
1 & 0 \cr
0 & \Sigma_{n-1}(f_{\alpha})}.
\]
It follows from Lemma \ref{lemma:matrixdifference} that the spectral
norm of $\tilde{ \Sigma}_n - \bar{\Sigma}_n$ converges to zero as
$n$ diverges.

\begin{figure}[b]

\includegraphics{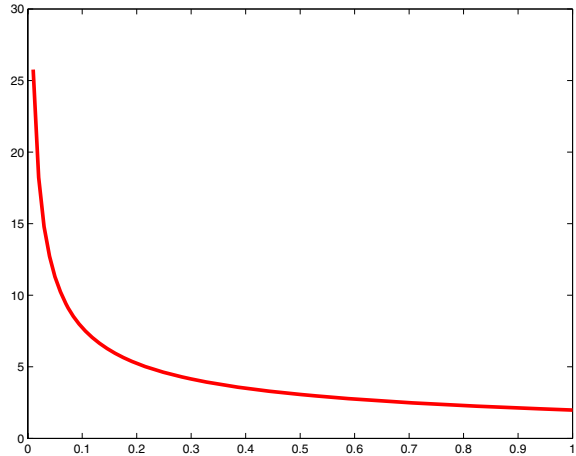}

\caption{Display of $C(f_{\alpha}, g_0)$. $x$-axis: $\alpha$. $y$-axis:
$C(f_{\alpha}, g_0)$.}
\label{fig:f0}
\end{figure}

Note that $\tilde{\mu} = \sqrt{a_n} \cdot\Sigma_{n-1}(g) \cdot\mu$
except for the first coordinate. Therefore, we expect model (\ref
{hj-model-1}) to be approximately equivalent to
\[
\tilde{X} = \sqrt{a_n} \cdot\Sigma_n(g_0) \cdot\mu+ \tilde{Z}
\qquad\mbox{where } \tilde{Z} \sim\mathrm{N}(0, \Sigma_{n}(f_{\alpha})).
\]
This is a special case of the cluster model we considered in Section
\ref{sec:cluster} with $f = f_{\alpha}$ and $g = g_0$, except that the
signal strength has been re-scaled by $\sqrt{a_n}$.
Therefore, if we calibrate the nonzero entries in $\mu$ as
%
%
\begin{equation} \label{hjsignal}
a_n^{-1/2} \cdot A_n = a_n^{-1/2} \cdot\sqrt{2 r \log n},
\end{equation}
then the detection boundary for the model is succinctly characterized by
\begin{eqnarray*}
r &=& \frac{1}{C(f_{\alpha}, g_0)} \cdot\rho^*(\beta),\\
C(f_{\alpha}, g_0) &=& \frac{1}{2\pi} \int_{-\pi}^{\pi} \frac{|g_0(\theta
)|^2}{f_{\alpha}(\theta)} \,d \theta= \frac{1}{\pi} \int_{-\pi
}^{\pi}
\frac{1 - \cos(\theta)}{f_{\alpha}(\theta)} \,d \theta.
\end{eqnarray*}
See Figure \ref{fig:f0} for the display of $C(f_{\alpha}, g_0)$. The
%
%
%
following theorem is proved in Section \ref{sec:proof}.
\begin{theorem} \label{thm:LB4}
Let $0 < \alpha_0 \leq\alpha< \frac{1}{2}$, $\beta\in(\frac
{1}{2},1)$, and $r \in(0,1)$. Assume $X$ is generated according to
model (\ref{hjmatrix}), with signal strength re-scaled as in (\ref
{hjsignal}). When $C(f_{\alpha}, g_0) \cdot r < \rho^*(\beta)$, the
null and alternative hypotheses merge asymptotically, and the sum of
types \textup{I} and \textup{II} errors of any test converges to $1$ as $n$ diverges
to infinity. When $ C(f_{\alpha}, g_0) \cdot r > \rho^*(\beta)$, if we
apply the iHC to $X$ with bandwidth $b_n = \log n$ and reject the
null when $\mathrm{iHC}_n^*(b_n, \Sigma_n) \geq(\log n)^2$, then the
type \textup{I} error converges to zero, and its power converges to 1.
\end{theorem}

\section{\texorpdfstring{Simulation study.}{Simulation study}} \label
{sec:Simul}

We conducted a small-scale empirical study to compare the performance
of iHC and standard HC. For iHC, we investigate two choices of
bandwidth: $b_n =1$ and $b_n = \log n$. In this section, we denote
standard HC, iHC with $b_n = 1$, and iHC with $b_n = \log n$ by HC,
HC-a and HC-b correspondingly.

The algorithm for generating data included the following four steps:
(1) Fix $n$, $\beta$, and $r$, let $m= n^{1 - \beta}$ and $A_n =
\sqrt
{2 r \log n}$. (2) Given a correlation matrix $\Sigma_n$, generate a
Gaussian vector $Z \sim\mathrm{N}(0, \Sigma_n)$. (3) Randomly draw $m$
integers $\ell_1 < \ell_2 < \cdots< \ell_m$ from $\{1, 2, \ldots,
n\}$
without replacement, and let $\mu$ be the $n$-vector such that $\mu_j
=A_n$ if $j\in\{\ell_1, \ell_2, \ldots, \ell_m\}$ and $0$ otherwise.
(4) Let $X = \mu+ Z$. Using data generated in this manner we explored
three parameter settings, (a)--(c), which we now describe.

In experiment (a), we took $n = 1000$ and $\Sigma_n(\rho)$ as the
tri-diagonal Toeplitz matrix generated by $f(\theta) = 1 + 2 \rho\cos
(\theta)$, $|\rho| < 1/2$. The
corresponding detection boundary was $r = \rho^*(\beta)/C(f)$ with
$C(f) = (2\pi)^{-1} \int_{-\pi}^{\pi}[1 - 2 \rho\cos(\theta
)]^{-1} \,d
\theta$. Consider all $\rho$ that range from $-0.45$ to $0.45$ with an increment
of $0.05$, and four pairs of parameters $(\beta, r) = (0.5, 0.2)$,
$(0.5, 0.25)$, $(0.55, 0.2)$ and $(0.55, 0.25)$. [Note that the
corresponding parameters $(m, A_n)$ are $(32,1.66)$, $(32,2.63)$,
$(22,1.66)$ and $(22,2.63)$]. For each triple $(\beta, r, \rho)$, we
generated data according to (1)--(4), applied HC, HC-a and HC-b to
both $Z$ and $X$ and repeated the whole process independently $500$ times.
As a result, for each triple $(\beta, r, \rho)$ and each procedure,
we got
$500$ HC scores that corresponded to the null hypothesis and $500$ HC
scores that corresponded to the alternative hypothesis.

\begin{figure}

\includegraphics{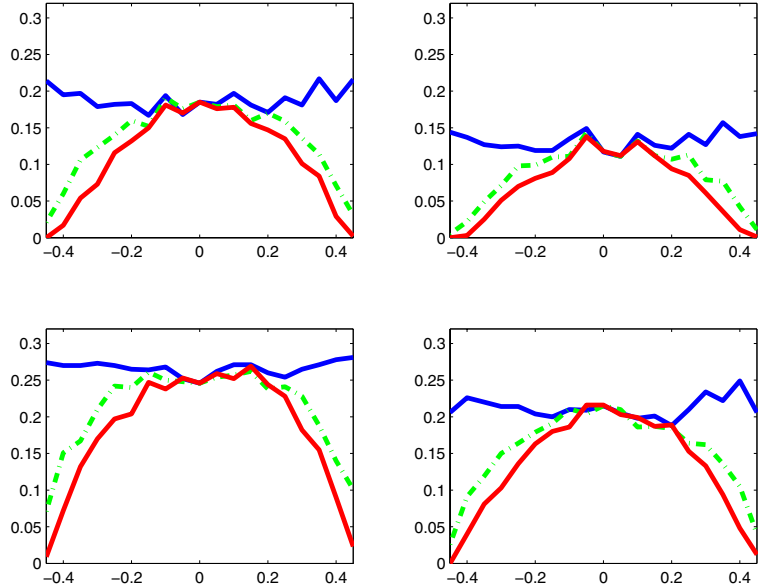}

\caption{Sum of types \textup{I} and \textup{II} errors as described in experiment
\textup{(a)}. From top to bottom then from left to right, $(\beta, r) =
(0.5,0.2)$, $(0.5,0.25)$, $(0.55,0.2)$, $(0.55,0.25)$. In each panel,
the $x$-axis displays $\rho$, and three curves (blue, dashed-green, and
red) display the sum of errors corresponding to HC, HC-a and HC-b.}
\label{fig:TriTotalerr}
\end{figure}

\begin{figure}

\includegraphics{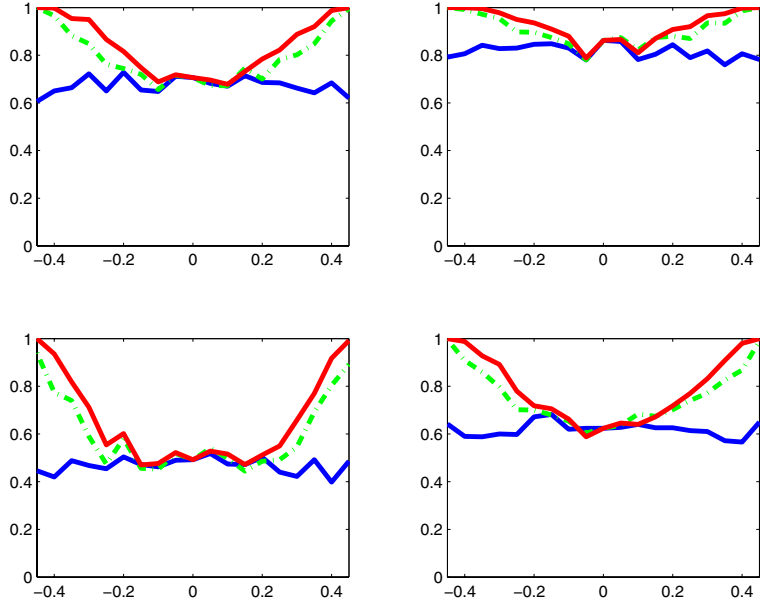}

\caption{Power as described in experiment \textup{(a)}. From top to bottom then
from left to right, $(\beta, r) = (0.5,0.2)$, $(0.5,0.25)$,
$(0.55,0.2)$, $(0.55,0.25)$. In each panel, the $x$-axis displays $\rho
$, and three curves (blue, dashed-green and red) display the power of
HC, HC-a and HC-b.}
\label{fig:TriPower}
\end{figure}

We report the results in two different ways. First,
we report the minimum sum of types I and II errors (i.e., the
minimum of the sum across all possible cut-off values) (see Figure \ref
{fig:TriTotalerr}). Second, we
pick the upper $10\%$ percentile of the $500$ HC scores corresponding
to the null hypothesis as a threshold (for later references, we call
this threshold the \textit{empirical threshold}) and calculate the
empirical power of the test (i.e., the fraction of HC scores
corresponding to the alternative hypothesis that exceeds the
threshold). The empirical thresholds are displayed in Table \ref
{table:cutoff} (to save space, only part of the thresholds are
reported), and the power is displayed in
Figure \ref{fig:TriPower}. Recall that in Theorem \ref{thm:HC} we
recommend $(\log n)^2$ as a cut-off point in the asymptotic setting.
For moderately large $n$, this cut-off point is conservative, and we
recommend the empirical threshold instead.

%
%

The results suggest that (1) iHC-b outperforms iHC-a, and iHC-a
outperforms HC. (2) As $|\rho|$ increases (note that a larger $|\rho|$
means a stronger correlation), the detection problem is increasingly easier,
and the advantage of iHC is increasingly prominent. (3) Under the null
hypothesis, the HC-b scores are usually smaller than those of HC and
HC-a. This is mainly due to the normalization term $\sqrt{2 b_n -1}$ in
the definition of iHC [see (\ref{DefineiHC})].

%
%
\begin{table}[b]
\caption{Display of empirical thresholds in experiment \textup{(a)} for
different $\rho$}
\label{table:cutoff}
\begin{tabular*}{\tablewidth}{@{\extracolsep{\fill}}lcccccccccc@{}}
\hline
$\bolds\rho$ & $\bolds-$\textbf{0.45} & $\bolds-$\textbf{0.35}
& $\bolds-$\textbf{0.25} & $\bolds-$\textbf{0.15} & $\bolds-$\textbf{0.05}
& \textbf{0.05} & \textbf{0.15} & \textbf{0.25} & \textbf{0.35}
& \textbf{0.45} \\
\hline
HC &3.059 & 2.851 & 2.913 & 2.892 &2.722 &2.835 & 2.742 &2.858 &2.834 &
3.032 \\
HC-a &2.919 & 2.858 & 2.924 & 2.837 &2.723 & 2.899& 2.713 & 2.826 &
2.677& 2.758 \\
HC-b &0.890 &0.847 &0.806 &0.773 & 0.769 & 0.775& 0.772 &0.761 &0.832
&0.859 \\
\hline
\end{tabular*}
\end{table}

%
%

%
%
\begin{figure}

\includegraphics{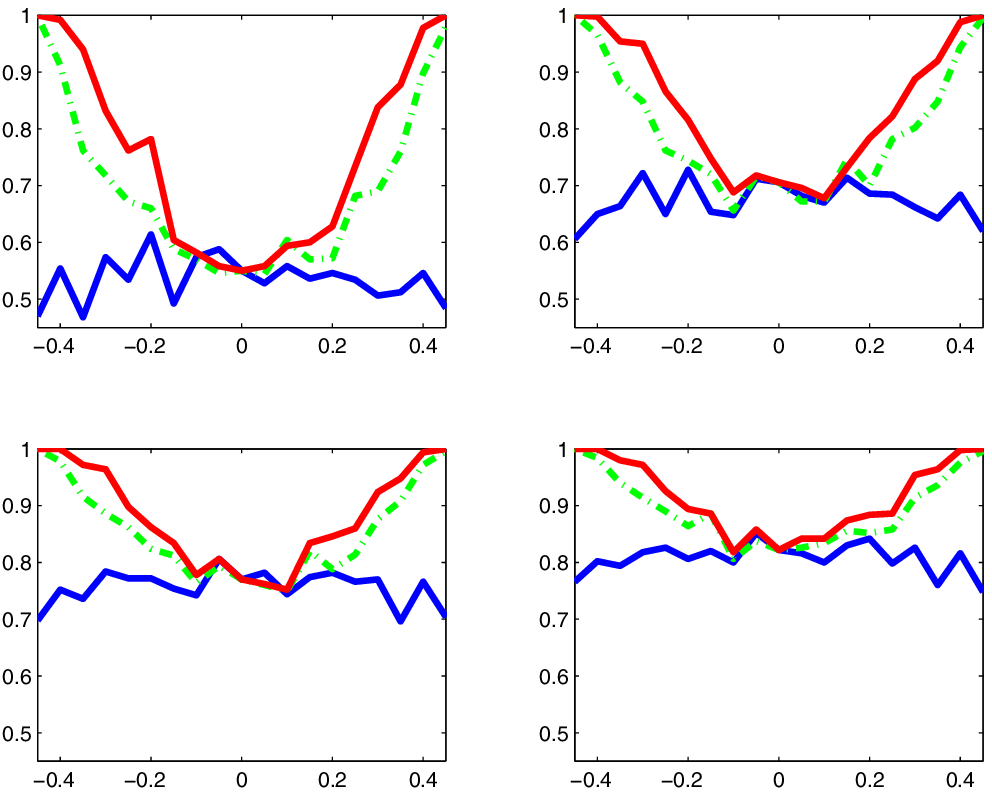}

\caption{Display of powers for different choices of cut-off value. Fix
$(\beta, r) = (0.5,0.2)$ as in experiment \textup{(a)}. From top to bottom
then from left to right, the cut-off values are the $5\%$, $10\%$, $15\%
$ and $20\%$ percentile of the $500$ HC scores corresponding to the
null hypothesis.
In each panel, the $x$-axis displays $\rho$, and three curves (blue,
dashed-green and red) display the power of HC, HC-a and HC-b. The
display suggest that, for different choices of cut-off value,
HC-b consistently outperforms HC-a, and HC-a
consistently outperforms HC.}
\label{fig:PowerAdd}
\end{figure}

We set the cut-off value as the $10\%$ percentile only for convenience.
Replacing $10\%$ by other percentage gives similar conclusion. See
Figure \ref{fig:PowerAdd} for details.

%
%
\begin{figure}

\includegraphics{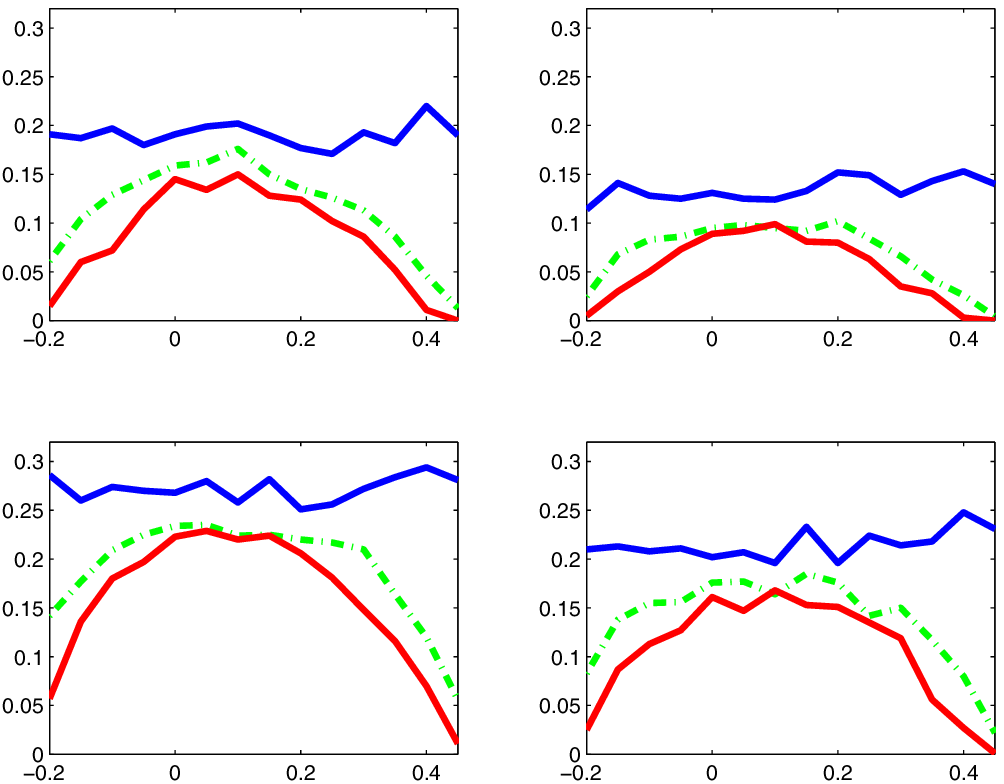}

\caption{Sum of types \textup{I} and \textup{II} errors as described in experiment
\textup{(b)}. From top to bottom then from left to right, $(\beta, r) =
(0.5,0.2)$, $(0.5,0.25)$, $(0.55,0.2)$, $(0.55,0.25)$. In each panel,
$x$-axis displays~$\rho$, and three curves (blue, dashed-green and red)
display the sum of errors corresponding to HC, HC-a and HC-b.}
\label{fig:5Totalerr}
\end{figure}

In experiment (b), we took $\Sigma_n$ to be the Toeplitz matrix
generated by
$f(\theta) = 1 + \frac{1}{2}\cos(\theta) + 2 \rho\cos(2 \theta)$ where
$\rho$ ranged from $-0.2$ to $0.45$ with an increment of $0.05$. (The
matrix $\Sigma_n$ is positive definite when $\rho$ is in this range.)
Other parameters are the same as in experiment (a). The minimum sums of
types I and II errors are reported in Figure \ref{fig:5Totalerr}.
The results suggest similarly that HC-b outperforms HC-a, and HC-a
outperforms HC.

In experiment (c), we investigated the behavior of HC-a/HC-b/HC for
larger $n$. We took
$(\beta, r) = (0.5, 0.25)$, $n = 500 \times(1, 2, 3, 4, 5)$ and
$\Sigma
_n$ as the tri-diagonal matrix in experiment (a) with $\rho= 0.4$.
The sum of types I and II errors is reported in Table \ref
{table:largern}. The results suggest that the performance of
HC-a/HC-h/HC improve when $n$ gets larger. (Investigation of the case
where $n$ was much larger than $2500$ needed much greater computer
memory, and so we omitted it.)

\section{\texorpdfstring{Discussion.}{Discussion}} \label{sec:Discu}

We have extended standard HC
to innovated HC by
building in the correlation structure.
The extreme diagonal entries of $\Sigma_n^{-1}$ play a key role in the
testing problem.
If the extreme value has finite upper and lower limits, $\bar{\gamma
}_0$ and $\underline{\gamma_0}$, then in the $\beta$--$r$ plane, the
detection boundary is bounded by the curves $r = \underline{\gamma
_0}^{-1} \cdot\rho^*(\beta)$ from above and $r = \bar{\gamma}_0^{-1}
\cdot\rho^*(\beta)$ from below.
When the correlation matrix is Toeplitz, the upper and lower limits
merge and equal the Wiener
interpolation rate $C(f)$. The detection boundary
is therefore $r=C(f)^{-1} \cdot\rho^*(\beta)$. The detection boundary
partitions the $\beta$--$r$ plane into a detectable region and an
undetectable region.
Innovated HC has asymptotically
full power for detection whenever $(\beta, r)$ falls into the interior
of the detectable region (we note, however, neither $\beta$ nor $r$ is
used to construct iHC). We call this the optimally adaptivity of
innovated higher criticism.

\subsection{\texorpdfstring{Connection to recent
literature.}{Connection to recent literature}}
The work complements that of\break Donoho and Jin \cite{DJ04} and Hall and
Jin \cite{HJ08}. The focus of \cite{DJ04} is standard HC and its
performance in the uncorrelated case. The focus of \cite{HJ08} is how
strong dependence may harm the effectiveness of standard HC; what could
be a remedy
was, however, not explored. The innovated HC proposed in the current
paper is optimal for both the model in \cite{DJ04} and that in \cite{HJ08}.

The work is related to that of Jager and Wellner \cite{Wellner} where
the authors proposed a family of
goodness-of-fit statistics for detecting sparse normal mixtures. The
work is also related to that of Meinshausen and Rice \cite{Rice} and of
Cai, Jin and Low~\cite{CJL}, where the authors focused on how to
estimate $\varepsilon_n$---the proportion of nonnull effects.

Recently, HC was also found to be useful for feature selection in
high-dimensional classification. See Donoho and Jin \cite{DJ08a,DJ08b},
Hall, Pittelkow and Ghosh \cite{HPG} and Jin \cite{JinPNAS}. The work
concerned the
situation where
there are relatively few samples containing a very large number of
features, out of which
only a small fraction is useful, and each useful feature contributes
weakly to the classification problem. In a related setting, Delaigle
and Hall \cite{DH} investigated HC for classification when the data is
non-Gaussian or dependent.

\subsection{\texorpdfstring{Future work.}{Future work}}\label{sec:further}
The work is also intimately connected to
recent literature on estimating covariance matrices. While the study is
focused more on situations where the correlation matrices can be
estimated using other approaches (e.g., \cite{Hongyu1,Goeman1,Goeman2}),
it can be generalized to cases
where the correlation matrix is unknown but can be estimated from data.
Cases where data on the covariance structure are available from other
time periods were discussed in Section \ref{estimatingsigma}, but even
if we stay within the confines of the current data, progress can be
made. In particular, it is noteworthy that it was shown in Bickel and
Levina \cite{Bickel} that when the correlation matrix has polynomial
off-diagonal decay, the matrix and its inverse can be estimated
accurately in terms of the spectral norm. In such situations we expect
the proposed approach to perform well once we combine it with that in
\cite{Bickel}.

%
\begin{table}[b]
\caption{Display of the sum of types \textup{I} and \textup{II} errors in experiment
\textup{(c)} for different $n$}
\label{table:largern}
\begin{tabular*}{\tablewidth}{@{\extracolsep{\fill}}lccccc@{}}
\hline
$\bolds n$ & \textbf{500} & \textbf{1000} & \textbf{1500} & \textbf{2000} & \textbf{2500} \\
\hline
HC & 0.201 & 0.144 & 0.115 & 0.123 & 0.098 \\
HC-a & 0.073 & 0.035 & 0.017 & 0.017 & 0.021 \\
HC-b & 0.033 & 0.005 & 0.004 & 0.003 & 0.002 \\
\hline
\end{tabular*}
\end{table}

Another interesting direction is to explore cases where the
correlation matrix does not have polynomial off-diagonal decay, but is sparse
in an unspecified pattern. This is a more challenging situation as
relatively little is known about the inverse of the correlation matrix.

Our study also opens opportunities for improving other recent
procedures. Take the aforementioned work on classification \cite
{DJ08a,DJ08b,HPG,JinPNAS}, for example. The approach derived in this
paper suggests ways of incorporating correlation structure into feature
selection, and therefore raises hopes for better classifiers. For
reasons of space, we leave explorations along these directions to
future study.

\section{\texorpdfstring{Proofs of main results.}{Proofs of main
results}} \label{sec:proof}

In this section we prove all theorems in preceding sections, except
Theorems \ref{thm:OHC} and \ref{thm:Toeplitz}. These two theorems are
the direct result of Theorems \ref{thm:LB} and \ref{thm:HC}, so we omit
the proofs. For simplicity, we drop the subscript $n$ whenever there is
no confusion.

\subsection{\texorpdfstring{Proof of (\protect\ref
{monotone}).}{Proof of (3.2)}} rewrite the second model
in (\ref{CompareExperiment}) as
$X + \xi= \mu+ \xi+ Z$,
where independently, $Z \sim\mathrm{N}(0, \Sigma)$,
$\xi\sim\mathrm{N}(0, \Delta)$, $\mu\sim G$ for $\Delta= \Sigma
_n^* -
\Sigma$ and some distribution $G$.
It suffices to show the monotonicity in the Hellinger distance. Denote
the density function of $\mathrm{N}(0, \Sigma)$ by $f(x) = f(x_1,
x_2, \ldots,
x_n)$, and write $dx_1 \, dx_2 \cdots dx_n$ as $dx$ for short. Then
the Hellinger distance corresponding to the second model in (\ref
{CompareExperiment}) can be written as
\[
h(\Sigma, \Delta, G) \equiv\int\sqrt{ E_{\Delta}\bigl( E_{G}
\bigl(f(x - \mu- \xi) \bigr) \bigr) \cdot E_{\Delta} \bigl(f(x - \xi
)\bigr) } \,dx.
\]
Note that by H\"older's inequality, $\sqrt{E[f] E[g]} \geq E[\sqrt
{fg}]$ for any positive and integrable functions $f$ and $g$. Using
Fubini's theorem, $h(\Sigma, \Delta, G)$ is not less than
\begin{eqnarray*}
&&\int\bigl[E_{\Delta} \sqrt{ \bigl( E_{G}f(x - \mu- \xi) \bigr) \cdot
\bigl( f(x - \xi) \bigr)}\bigr] \,dx \\
&&\qquad= E_{\Delta} \biggl[ \int\sqrt{
E_{G}f(x -
\mu- \xi) f(x - \xi) } \,dx \biggr].
\end{eqnarray*}
Note that $\int\sqrt{ E_{G}f(x - \mu- \xi) f(x - \xi) } \,dx
\equiv
\int\sqrt{ E_{G}f(x - \mu) f(x) } \,dx$
for any fixed $\xi$. It follows that
\begin{eqnarray*}
h(\Sigma, \Delta, G) &\geq& E_{\Delta} \biggl[ \int\sqrt{ E_{G}f(x -
\mu- \xi) f(x - \xi) } \,dx \biggr] \\
&=& \int\sqrt{ E_{G}f(x - \mu)
f(x) } \,dx,
\end{eqnarray*}
where the last term is the Hellinger distance corresponds to the first
model of (\ref{CompareExperiment}). Combining these results gives the
claim.

\subsection{\texorpdfstring{Proof of Theorem \protect\ref
{thm:LB}.}{Proof of Theorem 3.1}}

It is sufficient to show that the Hellinger distance between the joint
density of $X$ and $Z$ converges to zero as $n$ diverges to infinity.
By the assumption $\bar{\gamma}_0 r < \rho^*(\beta)$, we can
choose a
sufficiently small constant $\delta= \delta(r, \beta, \gamma_0)$ such
that $\bar{\gamma}_0 (1 - \delta)^{-2} r < \rho^*(\beta)$. Let
$\tilde
{\mu} = \mu/\sqrt{1 - \delta}$, let $U$ be the inverse of the Cholesky
factorization of $\Sigma$,
and let $\tilde{U}$ be the banded version of $U$,
\[
\tilde{U}(i,j) = \cases{
U(i,j), &\quad $|i - j| \leq\log^2(n)$, \cr
0, &\quad otherwise.}
\]
Model (\ref{Model0}) can be equivalently written as
%
%
\begin{equation} \label{Model2x0}
X = \tilde{\mu} + Z \qquad\mbox{where } Z \sim\mathrm{N}\bigl(0, (1 -
\delta)^{-1} \cdot\Sigma\bigr).
\end{equation}
The key to the proof is to compare model (\ref{Model2x0}) with the
following model:
%
%
\begin{equation} \label{Model2x}
X = \tilde{\mu} + Z \qquad\mbox{where } Z \sim\mathrm{N}(0, (\tilde
{U}' \tilde{U})^{-1}).
\end{equation}
In fact, by (\ref{monotone}), to establish the claim it suffices to
prove that (i) $ \tilde{U}' \tilde{U} \leq(1 - \delta)^{-1} \Sigma$
for sufficiently large $n$, and (ii) the Hellinger distance between the
joint density of $X$ and that of $Z$ associated with model (\ref
{Model2x}) tends to zero as $n$ diverges to infinity.

To prove the first claim, noting that $\Sigma= (U' U)^{-1}$, it
suffices to show $(1 - \delta) U' U \leq\tilde{U}' \tilde{U}$. Define
$W = U -\tilde{U}$ and observe that there is a generic constant $C > 0$
such that $\|\tilde{U}\| \leq C $ and $\|W\| \leq C$, whence $\|\tilde
{U}' \tilde{U} - \tilde{U}'
\tilde{U}\| = \|W' W + \tilde{U}' W + W' \tilde{U}\| \leq C \|W\|$.
Moreover, by \cite{Horn}, Theorem 5.6.9, for any symmetric matrix, the
spectral norm is no greater than the $\ell^1$-norm. In view of the
definitions of $W$ and $\Theta_n^*(\lambda, c_0, M)$, the $\ell^1$-norm
of $W$ is no greater than $(\log n)^{-2 (\lambda- 1)}$. Therefore, $\|
\tilde{U}' \tilde{U} - \tilde{U}'
\tilde{U}\| \leq C \|W\| \leq C (\log n)^{-2(\lambda-1)}$.
This, and the fact that all eigenvalues of $\tilde{U}' \tilde{U}$ are
bounded from below by a positive constant, imply the claim.

We now consider the second claim. Model (\ref{Model2x}) can be
equivalently written as $X = \tilde{U} \tilde{\mu} + Z$ where $Z
\sim
\mathrm{N}(0, I_n)$.
The key to the proof is that $\tilde{U}$ is a banded matrix and $\mu$
is a sparse vector where with probability converging to $1$, the
inter-distances of nonzero coordinates are no less than $3 (\log n)^2$
(see Lemma \ref{lemma:ell} for the proof). As a result the nonzero
coordinates of $\tilde{U} \tilde{\mu}$ are disjoint clusters of sizes
$O(\log^2 n)$ which simplifies the calculation of the Hellinger
distance. The derivation of the claim is summarized in Lemma \ref
{lemma:LB} which is stated and proved in the \hyperref[sec:appen]{Appendix}.

\subsection{\texorpdfstring{Proof of Theorem
\protect\ref{thm:HCAdd}.}{Proof of Theorem 4.1}}
Recall that $U_n$ is the function of $\Sigma_n$ defined by $U_n\Sigma_n
U_n'=I_n$. Put $Y = U_n X$, $\nu= U_n \mu$ and $Z = U_nz$. Model
(\ref
{ModelTransform1}) reduces to
%
%
\begin{equation} \label{ModelHCx1}
Y = \nu+ Z,\qquad Z \sim\mathrm{N}(0,I_n).
\end{equation}
Recalling that $\mathrm{HC}_n^*/\sqrt{2 \log\log n} \rightarrow1$ in
probability under $H_0$, it follows that $P_{H_0} \{ \mbox{Reject
$H_0$} \}$ tends to zero as $n$ diverges to infinity, and it suffices
to show $P_{H_1^{(n)}} \{\mbox{Accept $H_0$} \}
\rightarrow0$.\vspace*{1pt}

The key to the proof is to compare model (\ref{ModelHCx1}) with
%
%
\begin{equation} \label{ModelHCx2}
Y^* = \nu^* + Z \qquad\mbox{where } Z \sim\mathrm{N}(0,I_n),
\end{equation}
with $\nu^*$ having $m$ nonzero entries of equal strength $(1 - \delta
_n) A_n$ whose locations are randomly drawn from $\{1, 2, \ldots, n\}$
without replacement. By (\ref{ExplainUmu1}) and (\ref{ExplainUmu2}) and
the way $\tilde {\Theta }_n^*(\delta_n, b_n)$ is defined, we note that
$\nu_j \geq(1 - \delta _n) A_n$ for all $j \in\{ \ell_1, \ell_2,
\ldots, \ell_m\}$. Therefore,
%
%
\begin{equation} \label{ModelHCxx}
\mbox{signals in $\nu$ are both denser and stronger than those in
$\nu^*$}.
\end{equation}
Intuitively, standard HC applied to model (\ref{ModelHCx1}) is no
``less'' than that applied to model (\ref{ModelHCx2}).

We now establish this point. Let $\bar{F}_0(t)$ be the survival
function of
the central $\chi^2$-distribution $\chi^2_1(0)$, and
let $\bar{F}_n(t)$ and $\bar{F}_n^*$ be the empirical survival
function of $\{
Y_k^2\}_{k = 1}^n$ and $\{(Y_k^*)^2\}_{k =1}^n$, respectively. Using
arguments similar to those of Donoho and Jin \cite{DJ04} it can be
shown that standard HC applied to models (\ref{ModelHCx1}) and (\ref
{ModelHCx2}), denoted by $\mathrm{HC}_n^{(1)}$ and $\mathrm{HC}_n^{(2)}$
for short, can be rewritten as
\begin{eqnarray*}
\mathrm{HC}_n^{(1)} &=& \sup_{t\dvtx1/n \leq\bar{F}_0(t) \leq1/2}
\biggl\{ \frac{\sqrt{n} (\bar{F}_n(t) - \bar{F}_0(t))}{\sqrt{\bar
{F}_0(t) F_0(t)}}
\biggr\},
\\
\mathrm{HC}_n^{(2)} &=& \sup_{t\dvtx1/n \leq\bar{F}_0(t) \leq1/2}
\biggl\{
\frac{\sqrt{n} (\bar{F}_n^*(t) - \bar{F}_0(t))}{\sqrt{\bar
{F}_0(t) F_0(t)}}
\biggr\},
\end{eqnarray*}
respectively.\vspace*{1pt}
The key fact is now that the family of noncentral $\chi
^2$-distribution $\{\chi_1^2(\delta), \delta\geq0\}$ is a monotone
likelihood ratio family (MLR),
that is, for any fixed $x$ and $\delta_2 \geq\delta_1 \geq0$,
$P\{\chi_1^2(\delta_2) \geq x\} \geq P\{\chi^2_1(\delta_1) \geq x\}$.
Consequently, it follows from (\ref{ModelHCxx}) and mathematical
induction that for any $x$ and $t$,
$P\{\bar{F}_n^*(t) \geq x\} \geq P\{\bar{F}_n(t) \geq x\}$.
Therefore, for any
fixed $x > 0$,
%
%
\begin{equation} \label{ModelHCxxx}
P \bigl\{ \mathrm{HC}_n^{(1)} < x \bigr\} \leq P \bigl\{ \mathrm{HC}_n^{(2)} <
x \bigr\}.
\end{equation}
Finally, by an argument similar to that of Donoho and Jin \cite{DJ04},
Section 5.1, the second term in (\ref{ModelHCxxx}) with $x = (1+a)
\sqrt
{2 \log\log n}$ tends to zero as $n$ diverges to infinity. This
implies the claim.

\subsection{\texorpdfstring{Proof of Theorem \protect\ref
{thm:HC}.}{Proof of Theorem 4.2}}
In view of Lemma \ref{lemma:HCNull}, it suffices to show that
$P_{H_1^{(n)}} \{\mbox{Accept $H_0$}\} \rightarrow0$. Put\vspace*{1pt} $\bar{U} =
\bar
{U}(b_n)$, $V = V_n(b_n)$, $Y = V X$, $\nu= V \mu$, $\tilde{Z} = V
Z$. Model (\ref{ModelTransformx}) reduces to
%
%
\begin{equation} \label{HCmodelx1}
Y = \nu+ \tilde{Z} \qquad\mbox{where } \tilde{Z} \sim\mathrm{N}(0,
\bar{U}' \bar{U}).
\end{equation}
Let $\bar{F}_n(t)$ and $\bar{F}_0(t)$ be the empirical survival
function of
$\{
Y_k^2\}_{k = 1}^n$ and the survival function of $\chi^2_1(0)$,
respectively. Let $q = q(\beta, r) =\break \min\{ (\beta+ \bar{\gamma}_0
r)^2/(4 \bar{\gamma}_0 r), 4 \bar{\gamma}_0 r\}$ and set $t_n^* =
\sqrt
{2 q \log n}$. Since $\bar{\gamma}_0 r < \rho^*(\beta)$, then it
can be
shown that $0 < q < 1$ and $n^{-1} \leq\bar{F}_0(t_n^*) \leq1/2$ for
sufficiently large $n$. Using an argument similar to that in the proof
of Theorem \ref{thm:HCAdd},
\begin{eqnarray*}
\mathrm{iHC}_n^* &=& \sup_{s\dvtx1/n \leq\bar{F}_0(s) \leq1/2}
\frac{
\sqrt
{n}(\bar{F}_n(s) - \bar{F}_0(s))}{\sqrt{(2 b_n-1)\bar{F}_0(s) (1 -
\bar{F}_0(s))}}\\
&\geq& \frac{\sqrt{n}(\bar{F}_n(t_n^*) - \bar{F}_0(t_n^*))}{\sqrt{(2
b_n -1)
\bar{F}
_0(t_n^*) (1 - \bar{F}_0(t_n^*))}}
\end{eqnarray*}
and it follows that
%
%
\begin{equation} \label{HCmodelxx}
P\{ \mathrm{iHC}_n^* \leq\log^{3/2}(n)\} \leq P\biggl\{ \frac{\sqrt
{n}(\bar{F}
_n(t_n^*) - \bar{F}_0(t_n^*))}{\sqrt{(2 b_n -1) \bar{F}_0(t_n^*) (1
- \bar{F}
_0(t_n^*))}} \leq\log^{3/2}(n)\biggr\}.\hspace*{-37pt}
\end{equation}
It remains to show that the right-hand side of (\ref{HCmodelxx}) is
algebraically small. The proof needs detailed calculations summarized
in Lemma \ref{lemma:HCproof} which is stated and proved in
the \hyperref[sec:appen]{Appendix}.

\subsection{\texorpdfstring{Proof of Theorem
\protect\ref{thm:cluster}.}{Proof of Theorem 6.1}}
Inspection of the proof of Theorems \ref{thm:LB} and~\ref{thm:HC}
reveals that the condition that $\Sigma_n$ is a correlation matrix and
that $\Sigma_n \in\Theta_n^*(\lambda, c_0$, $M)$ in those theorems can
be relaxed. In particular, $\Sigma_n$ need not have equal diagonal
entries, and the decay condition on $\Sigma_n$ can be replaced by a
weaker condition that concerns the decay of $U_n$ (the inverse of the
Cholesky factorization of $\Sigma_n$), specifically
\[
|U_n(i,j)| \leq M (1 + |i - j|^{\lambda})^{-1}.
\]

Let $U_n(f)$ be the inverse of the Cholesky factorization of $\Sigma
_n(f)$, and define $\tilde{U}_n = U_n(f) \Sigma_n(g)$. Since $\Sigma
_n(g)$ is a lower triangular matrix with positive diagonal entries,
then it is seen that $\tilde{U}_n$ is the inverse of the Cholesky
factorization of $\tilde{\Sigma}_n$. By Lemma \ref{lemma:Udecay},
$U_n(f)$ has polynomial off-diagonal decay with
the parameter $\lambda$. It follows that $\tilde{U}_n$ decays at the
same rate. Applying Theorems \ref{thm:LB} and \ref{thm:HC}, we see that
all that remains to prove is that
%
%
\begin{equation} \label{clustercheck}
{\max_{\sqrt{n} \leq k \leq n - \sqrt{n}}}|\tilde{\Sigma
}_n^{-1}(k,k) -
C(f, g)|\rightarrow0.
\end{equation}

By \cite{Bottcher}, Theorem 2.15, for any $\sqrt{n} \leq k \leq n -
\sqrt{n}$, $k - K \leq j \leq k + K$ and $1 \leq\lambda' < \lambda$,
\[
\bigl|\Sigma_n^{-1}(f)(k, j) - \bigl(\Sigma_n(1/f)\bigr)(k,j)\bigr| = o\bigl(n^{-(1 - \lambda')/2}\bigr).
\]
Since $\tilde{\Sigma}_n^{-1} = \Sigma_n(\bar{g}) \cdot\Sigma_n^{-1}(f)
\cdot\Sigma_n(g)$,
it follows that
$\sup_{\sqrt{n} \leq k \leq n - \sqrt{n}} |\tilde{\Sigma
}_n^{-1}(k,k) -
(\Sigma_n(\bar{g}) \cdot\Sigma_n(1/f) \cdot\Sigma_n(g))(k,k)|
\rightarrow
0$. Moreover, direct calculations show that $(\Sigma_n(\bar{g}) \cdot
\Sigma_n(1/f) \cdot\Sigma_n(g))(k,k) = C(f,g)$, $\sqrt{n} \leq k
\leq
n - \sqrt{n}$. Combining these results gives (\ref{clustercheck}) and
concludes the proofs.

\subsection{\texorpdfstring{Proof of Theorem \protect\ref
{thm:LB4}.}{Proof of Theorem 7.1}}
Consider the first claim. It suffices to show that the Hellinger
distance between $\tilde{X}$ and $\tilde{Z}$ in model (\ref
{hj-model-1}) tends to zero as $n$ diverges to infinity. Since $
C(f_{\alpha}, g_0) \cdot r < \rho^*(\beta)$, there is a small constant
$\delta> 0$ such that $ (1 - \delta)^{-1} \cdot C(f_{\alpha}, g_0)
\cdot r < \rho^*(\beta)$. Using Lemma \ref{lemma:f0}, we see that
$\Sigma_{n-1}(f_{\alpha})$ is a positive matrix the smallest eigenvalue
of which is bounded away from zero. It follows from Lemma \ref
{lemma:matrixdifference} and basic algebra that
$\tilde{\Sigma} \geq(1 - \delta) \bar{\Sigma}_n$ for sufficiently
large $n$.
Compare model (\ref{hj-model-1}) with
%
%
\begin{equation}\label{hj-model-2}
X^{*} = \tilde{\mu} + Z^{*} \qquad\mbox{where } Z^{*} \sim\mathrm{N}
\bigl(0, (1 - \delta) \bar{\Sigma}\bigr).
\end{equation}
By the monotonicity of Hellinger distance at (\ref{monotone}), it
suffices to show
that the Hellinger distance between $X^*$ and $Z^*$ tends to zero as
$n$ diverges to infinity.

Now, by the definition of $\tilde{\mu}$,
$\tilde{\mu} - \sqrt{a_n} \cdot\Sigma_{n}(g_0) \cdot\mu= (\mu_n,
\sqrt{a_n} \cdot\mu_n, 0, \ldots, 0)'$.
Since $P\{\mu_n \neq0\} = o(1)$
then, except for an event with negligible probability, \mbox{$\tilde{\mu}
=\bar{\mu}$}.
Therefore, replacing $\tilde{\mu}$ by $\sqrt{a_n} \cdot\Sigma_{n}(g_0)
\cdot\mu$ in model (\ref{hj-model-2})
alters the Hellinger distance only negligibly.
Note that the first coordinate of $X^*$ is uncorrelated with all other
coordinates, and its mean equals zero with probability converging to
$1$, so removing it from the model only has a negligible effect on the
Hellinger distance. Combining these properties, model (\ref
{hj-model-2}) reduces to the following with only a negligible
difference in the Hellinger distance:
\begin{eqnarray*}
X^{*}(2\dvtx n) &=& \Sigma_{n-1}(g_0) \bigl(\sqrt{a_n} \cdot\mu(2\dvtx n)\bigr) +
Z^{*}(2\dvtx n),\\
Z^*(2\dvtx n) &\sim& \mathrm{N}\bigl(0, (1 - \delta) \Sigma
_{n-1}(f_{\alpha})\bigr),
\end{eqnarray*}
where $X(2\dvtx n)$ denotes the vector $X$ with the first entry removed.
Dividing both sides by $\sqrt{1 - \delta}$, this reduces to the
following model:
%
%
\begin{eqnarray}\label{hj-model-3}
\tilde{X}(2\dvtx n) &=& \Sigma_{n-1}(g_0) \frac{\sqrt{a_n} \cdot\mu
(2\dvtx n)}{\sqrt{1 - \delta}} + \tilde{Z}(2\dvtx
n),\nonumber\\[-8pt]\\[-8pt]
\tilde{Z}(2\dvtx n)
&\sim&
\mathrm{N}(0, \Sigma_{n-1}(f_{\alpha})),\nonumber
\end{eqnarray}
which is in fact model (\ref{clustermodel}) considered in Section \ref
{sec:cluster}.
It follows from (\ref{hjsignal}) that $\sqrt{a_n} \cdot\mu
(2\dvtx n)/\sqrt
{1 - \delta}$ has $m$ nonzero coordinates each of which equals $\sqrt{2
(1 - \delta)^{-1} r \log n}$. Comparing model (\ref{hj-model-3}) with
model (\ref{clustermodel}) and recalling that $(1 - \delta)^{-1}
\cdot
r \cdot C(f_{\alpha}, g_0) < \rho^*(\beta)$, the claim follows from
Theorem \ref{thm:cluster}.

Consider the second claim. Since $ C(f_{\alpha},g_0) \cdot r > \rho
^*(\beta)$, then there is a small constant $\delta> 0$ such that $(1 -
\delta) \cdot r \cdot C(f_{\alpha}, g_0) > \rho^*(\beta)$. Let $U_n$
be the inverse of the Cholesky factorization of
$\Sigma_n$, and let $\bar{U}_n(b_n)$ and $V_n(b_n)$ be as defined right
below (\ref{DefineUbar}).
Write model (\ref{hjmatrix}) equivalently as
\[
V X = V \mu+ VZ \qquad\mbox{where $VZ \sim\mathrm{N}(0, \bar{U}'(b_n)
\bar
{U}(b_n))$}.
\]
Recall that $\bar{U}'(b_n) \bar{U}(b_n)$ is a banded correlation matrix
with bandwidth $2b_n -1$. Let $\ell_1, \ell_2, \ldots, \ell_m$ be the
$m$ locations of nonzero means of $\mu$. By an argument similar to that
in the proof of Theorem \ref{thm:HC}, all remains to show is that,
except for an event with negligible probability,
%
%
\begin{eqnarray} \label{hjv}
(V \mu)_{k} \geq\sqrt{2 r' \log n} \hspace*{85pt}\nonumber\\[-8pt]\\[-8pt]
\eqntext{\mbox{for some constant $r'
> \rho^*(\beta)$ and all $k \in\{\ell_1, \ell_2, \ldots, \ell_k\}$}.}
\end{eqnarray}

We now show (\ref{hjv}).
First, by Lemma \ref{lemma:Transformmean} and (\ref{hjsignal}), except
for an event with negligible probability,
\[
(V \mu)_{k} \geq(1 - \delta)^{1/4} \cdot\bigl(a_n \cdot\Sigma
_n(k,k)\bigr)^{-1/2} \cdot A_n,\qquad k \in\{\ell_1, \ell_2, \ldots,
\ell
_m\}.
\]
Second, by the way $\tilde{\Sigma}_n$ is defined,
\[
(a_n \Sigma_n^{-1})(k,k) = \bigl(\Sigma_{n}(g_0) \cdot\tilde{\Sigma
}_n^{-1} \cdot\Sigma_n(\bar{g}_0)\bigr)(k,k)\qquad \mbox{for all $k
\geq2$},
\]
and by the way $\bar{\Sigma}_n$ is defined and Lemma \ref
{lemma:matrixdifference}, for sufficiently large $n$,
\[
\tilde{\Sigma}_n^{-1} \geq(1 - \delta)^{-1/2} \bar{\Sigma
}_n^{-1},
\]
and so
\[
\Sigma_{n}(g_0) \tilde{\Sigma}_n^{-1}
\Sigma
_n(\bar{g}_0) \geq(1 - \delta)^{1/2} \Sigma_{n}(g) \bar{\Sigma
}_n^{-1} \Sigma_n(\bar{g}).
\]
Last, by \cite{Bottcher}, Theorem 2.15, $|(\Sigma_n(g_0) \cdot\bar
{\Sigma}_n^{-1} \cdot\Sigma_n(\bar{g}_0))(k,k) - C(f_{\alpha},
g_0)| =
o(1)$ when $\min\{k, n-k\}$ is sufficiently large. Combining these
results gives (\ref{hjv}) with $r' = (1 - \delta) \cdot r \cdot
C(f_{\alpha}, g_0)$, and the claim follows directly.

\begin{appendix}
\section*{Appendix} \label{sec:appen}

\subsection{\texorpdfstring{Statement and proof of Lemma
\protect\ref{lemma:Udecay}.}{Statement and proof of Lemma
A.1}}

\begin{lemma} \label{lemma:Udecay}
Fix $\lambda> 1$, $c_0 > 0$, and $M > 0$. For any sequence of matrices
$\Sigma_n$, $n\geq1$, such that $\Sigma_n \in\Theta_n^*(\lambda, c_0,
M)$, let $U_n$ be the inverse of the Cholesky factorization of $\Sigma
_n$. Then there is a constant $C = C(\lambda, c_0, M) > 0$ such that,
for any $n$ and any $1 \leq j, k \leq n$,
\[
|\Sigma_n^{-1}(j,k)| \leq C \cdot(1 + |j - k|)^{-\lambda},\qquad
|U_n(j,k)| \leq C \cdot(1 + |j - k|)^{-\lambda}.
\]
\end{lemma}
\begin{pf}
When $\lambda= 1$, the first inequality continues to hold, and the
second holds if we adjoin a $\log n$ factor to the right-hand side.

As a prelude to giving the proof we state the following result, taken
directly from \cite{Sun1}.
Let $\mathbb{Z}$ be the set of all integers. Write $\ell^2$ for the set
of summable sequences $x = \{x_k\}_{k \in\mathbb{Z}}$, and let $A =
(A(j,k))_{j,k \in\mathbb{Z}}$ be an infinite matrix. Also, let $|x|_2$
be the $\ell^2$-vector norm of $x$, and $\|A\|$ be the operation norm
of $A$: $\|A\|={\sup_{x \dvtx| x |_2 = 1}}|A x|_2$. Fixing positive
constants $\lambda$, $M$ and $c_0$, we define
the class of matrices
%
%
\setcounter{equation}{0}
\begin{eqnarray} \label{DefineQ}
&&\Theta_{\infty}(\lambda, c_0, M) = \biggl\{A = (A(j,k))_{j,k \in
\mathbb
{Z}}\dvtx|A(j,k)| \nonumber\\[-8pt]\\[-8pt]
&&\hspace*{92.37pt}\leq\frac{M}{(1 + |j - k|)^{\lambda}}, \|A\|
\geq
c_0 \biggr\}.\nonumber
\end{eqnarray}
\begin{lemma} \label{lemma:matrixdecay}
Fix $\lambda> 1$, $c_0 > 0$, and $M > 0$. For any matrix $A \in\Theta
_{\infty}(\lambda, M)$, there is a constant $C > 0$, depending only on
$\lambda$, $M$ and $c_0$, such that $|A^{-1}(j,k)| \leq C\cdot(1 + |j
- k|)^{-\lambda}$.
\end{lemma}

Next we consider the first claim in Lemma \ref{lemma:Udecay}. Construct
an infinite matrix $\Sigma_{\infty}$ by arranging the finite matrices
along the diagonal, and note that the inverse of $\Sigma_{\infty}$ is
the matrix
formed by arranging the inverse of the finite matrices along the
diagonal. Since $\Sigma_{\infty}(i,j) \leq M (1 + |i - j|^{\lambda
})^{-1}$, then applying Lemma \ref{lemma:matrixdecay} gives the claim.

Consider the second claim. It suffices to show that $|U_n(k,j)| \leq C
/(1 + |k - j|^{\lambda})$ for all $1 \leq j < k \leq n$.
Denote the first $k \times k$ main diagonal sub-matrix of $\Sigma_{n}$
by $\Sigma_{(k)}$, the $k$th row of $\Sigma_{(k)}$ by $(\xi_{k-1}',
1)$, and the $k$th row of $U_n$ by $u_k'$. It follows from direct
calculations that
%
%
\begin{equation} \label{Udecay1}
u_k' = \bigl( 1 - \xi_{k-1}' \Sigma_{(k-1)}^{-1} \xi_{k-1}
\bigr)^{-1/2} \cdot\bigl(\xi_{k-1}' \Sigma_{(k-1)}^{-1}, 1\bigr).
\end{equation}
At the same time, by (\ref{Udecay1}) and basic algebra,
%
%
\begin{equation} \label{Udecay2}
\bigl(1 - \xi_{k-1}' \Sigma_{(k)}^{-1} \xi_{k-1} \bigr)^{-1} \leq u_k'
u_k = \Sigma_{(k)}^{-1}(k,k).
\end{equation}
Combining (\ref{Udecay1}) and (\ref{Udecay2}) gives
%
%
\begin{equation} \label{Udecay3}
|U_n(k,j)| = |u_{k}(j)| \leq C \bigl| \bigl(\Sigma_{(k-1)}^{-1}
\xi_{k-1}\bigr)_j\bigr|,\qquad
1 \leq j \leq k -1.
\end{equation}
Now, by Lemma \ref{lemma:matrixdecay}, $|\Sigma_{(k-1)}^{-1}(j,s)|
\leq
C(1 + |j - s|^{\lambda})^{-1}$ for all $1 \leq i, j \leq k-1$. Note
that $|\xi_{k-1}(s)| \leq C (1 + |s - k|^{\lambda})^{-1}$, $1 \leq s
\leq n$ and $\lambda> 1$. It follows from basic algebra that
%
%
\begin{equation} \label{Udecay4}
\bigl|\bigl(\Sigma_{(k-1)}^{-1} \xi_{k-1}\bigr)_j\bigr| \leq\sum_{s = 1}^n \frac{C}{(1 +
|j - s|^{\lambda})(1 + |s - k|^{\lambda})} \leq\frac{C}{1 + |k -
j|^{\lambda}}.
\end{equation}
Inserting (\ref{Udecay4}) into (\ref{Udecay3}) gives the claim.
\end{pf}

\subsection{\texorpdfstring{Statement and proof of Lemma
\protect\ref{lemma:Transformmean}.}{Statement and proof of Lemma
A.3}}

\begin{lemma} \label{lemma:Transformmean}
Fix $c_0 > 0$, $\lambda\geq1$, and $M > 0$. Consider a sequence of
bandwidths $b_n$ that tends to infinity. Let $\{\ell_1, \ell_2,
\ldots,
\ell_m\}$ be the $m$ random locations of signals in $\mu$, arranged in
the ascending order. For sufficiently large $n$, there is a constant $C
= C(c_0, \lambda, M)$ such that, except for an event with
asymptotically vanishing probability,
\[
(V_n(b_n) \mu)_k \geq\bigl(1 - C b_n^{1/2 - \lambda} + o(1)\bigr) \cdot\sqrt
{\Sigma_n^{-1}(k,k)} \cdot A_n\qquad \forall k \in\{\ell_1,
\ell
_2, \ldots, \ell_m\},
\]
for all $\Sigma_n \in\Theta_n^*(\lambda, c_0, M)$, where $o(1)$ tends
to zero algebraically fast.
\end{lemma}
\begin{pf}
To derive the lemma, note that we may assume without loss of generality
that $\ell_1 < \ell_2 < \cdots< \ell_m$. By Lemma \ref{lemma:ell},
except for an event with negligible probability, $\ell_1 \geq b_n$,
$\ell_m \leq n -b_n$, and the inter-$\ell_j$ distances are not less
than $C \log n \cdot n^{2 \beta- 1}$. For any $k \in\{\ell_1, \ell_2,
\ldots, \ell_m\}$, let
$d_k = (\sum_{j = k}^{k + b_n -1} u_{jk}^2 )^{-1/2}$.
By the way $\bar{U}(b_n)$ is defined,
%
%
\begin{eqnarray} \label{Transformmean0}
(\bar{U}'(b_n) U \mu)_k &=& d_k \sum_{s, j =1}^n \tilde{u}_{k s} u_{s j}
\mu_j \nonumber\\[-8pt]\\[-8pt]
&=& d_k \Biggl[ \sum_{s, j = 1}^n u_{ks} u_{sj} \mu_j - \sum
_{s,j=1}^n (u_{ks} - \tilde{u}_{ks}) u_{sj} \mu_j \Biggr].\nonumber
\end{eqnarray}

Consider $d_k$ first. Write
\[
1/d_k^2 = \sum_{j = k}^{k + b_n -1} u_{jk}^2 = \sum_{j=k}^{n} u_{jk}^2
- \sum_{j = k-b_n}^n u_{jk}^2.
\]
First, $U'U = \Sigma^{-1}$, $\sum_{j=k}^{n} u_{jk}^2 = (U'U)(k,k) =
(\Sigma^{-1})(k,k)$. Second, by the polynomial off-diagonal decay of
$U$ and basic calculus,
\[
\sum_{j = k + b_n}^n u_{jk}^2 \leq C \sum_{j = k+b_n}^n \frac{1}{1 + |j
- k|^{\lambda}} \leq C b_n^{1 - 2 \lambda}.
\]
Last, note that the quantities $\Sigma^{-1}(k,k)$ are uniformly bounded
away from zero and infinity. Combining these results gives
%
%
\begin{equation} \label{Transformmean1}
\bigl|d_k - \sqrt{ \Sigma^{-1}(k,k)}\bigr| \leq C b_n^{1 - 2 \lambda}.
\end{equation}

Consider $\sum_{s, j = 1}^n u_{ks} u_{sj} \mu_j$ next.
Recall that $\mu_j = A_n$ when $j \in\{\ell_1, \ell_2, \ldots,
\ell_m\}
$ and $\mu_j = 0$ otherwise. Since $U'U = \Sigma^{-1}$,
\[
\sum_{s, j = 1}^n u_{ks} u_{sj} \mu_j = \sum_{j = 1}^n (\Sigma
^{-1})(k,j) \mu_j = A_n \Sigma^{-1}(k,k) + A_n \sum_{\ell_s \neq k}
\Sigma^{-1}(k, \ell_s).
\]
Define $L_n = n^{\beta- 1/2}$.
By Lemma \ref{lemma:ell}, except for an event with negligible
probability, the inter-distance of $\ell_j$ is no less than $L_n$. So
by the polynomial off-diagonal decay of $\Sigma^{-1}$, the second term
is algebraically small. Therefore,
%
%
\begin{equation} \label{Transformmean2}
\sum_{s, j = 1}^n u_{ks} u_{sj} \mu_j = A_n [(\Sigma^{-1})(k,k) +
o(b_n^{1 - \lambda})].
\end{equation}

Last, we consider $\sum_{s,j=1}^n (u_{ks} - \tilde{u}_{ks}) u_{sj}
\mu_j$.
Direct calculations show that
\[
\bigl|\bigl((U - \tilde{U})' U\bigr)(k, j)\bigr|
\leq
\cases{
\dfrac{C}{1 + |k - j|^{\lambda}}, &\quad $\lambda> 1$, \cr
\dfrac{C \log n}{1 + |k - j|^{\lambda}}, &\quad $\lambda= 1$,}
\]
so by a similar argument,
\begin{eqnarray*}
\Biggl|\sum_{s,j=1}^n (u_{ks} - \tilde{u}_{ks}) u_{sj} \mu_j
\Biggr| &=&
\Biggl|\sum_{j = 1}^n \bigl((U - \tilde{U})' U\bigr)(k,j) \mu_j \Biggr| \\
&\leq& A_n
\cdot\bigl((U - \tilde{U})' U\bigr)(k,k) + o(1),
\end{eqnarray*}
where $o(1)$ is algebraically small.
Moreover, by the inequality,
\[
\bigl((U - \tilde{U})' U\bigr)(k,k) \leq\sum_{s = 1}^n |(u_{ks} - \tilde
{u}_{ks}) u_{sk}| \leq b_n^{ 1/2 - \lambda}
\]
and the claim follows.
\end{pf}

\subsection{\texorpdfstring{Statement and proof of Lemma \protect
\ref
{lemma:empiricalprocess}.}{Statement and proof of Lemma A.4}} \label{sub:empprocess}

Let $p_1,\ldots,p_N$ be $N$ independent and identically distributed
data from $U(0,1)$,
and $F_N(t)$ be the empirical cdf. The normalized uniform stochastic
process is defined as
\[
\mathbb{W}_N(t) = \sqrt{N}[F_N(t) - t]/\sqrt{t(1 - t)}.
\]
\begin{lemma} \label{lemma:empiricalprocess}
There is a generic constant $C > 0$ such that for sufficiently large~$n$,
\[
P \Bigl\{ {\sup_{1/n \leq t \leq1/2}} |\mathbb{W}_N(t)| \geq C (\log
n)^{3/2} \Bigr\} \leq C n^{-C}.
\]
\end{lemma}
\begin{pf}
To derive this result, note that by the Hungarian construction \cite
{Csorgo}, there is a Brownian bridge $\mathbb{B}(t)$ such that
\[
P \biggl\{ \sup_{1/n \leq t \leq1/2} \bigl|\sqrt{N}\bigl(F_N(t) - t\bigr) - \mathbb
{B}(t)\bigr| \geq\frac{ C (\log N + x)}{\sqrt{N}} \biggr\} \leq C e^{-C
x},
\]
where $C > 0$ are generic constants. Noting that $1/\sqrt{t(1 - t)}
\leq\sqrt{n} \leq C \sqrt{N \log N}$ when $1/n \leq t \leq1/2$, it
follows that
%
%
\begin{eqnarray} \label{empiricalprocess1}
&&P \biggl\{ \sup_{1/n \leq t \leq1/2} \biggl|\frac{\sqrt{N}(F_N(t) -
t) - \mathbb{B}(t)}{\sqrt{t(1 - t)}} \biggr| \geq C (\log N)^{1/2}
(\log N + x) \biggr\}\nonumber\\[-8pt]\\[-8pt]
&&\qquad \leq C e^{-C x}.\nonumber
\end{eqnarray}
At the same time, by \cite{Wellnerbook}, page 446,
%
%
\begin{equation} \label{empiricalprocess2}
P \biggl\{ \sup_{1/n \leq t \leq1/2} \biggl|\frac{ \mathbb
{B}(t)}{\sqrt{t(1 - t)}} \biggr| \geq C (\log N)^{1/2} x \biggr\}
\leq
C\log N \cdot e^{-C x}.
\end{equation}
Combining (\ref{empiricalprocess1}) and (\ref{empiricalprocess2}),
taking $x = C \log N$ and using the triangle inequality, we deduce the
lemma.
\end{pf}

\subsection{\texorpdfstring{Statement and proof of Lemma
\protect\ref{lemma:HCNull}.}{Statement and proof of Lemma A.5}}

\begin{lemma} \label{lemma:HCNull}
Take the bandwidth to be $b_n = \log n$ and suppose $H_0$ is true.
Then, except for an algebraically small probability,
$\mathrm{iHC}_n^*(b_n) \leq C (\log n)^{3/2}$ for some constant $C > 0$,
uniformly for all correlation matrices.
\end{lemma}
\begin{pf}
To derive the lemma, note that we may assume without loss of gen\-erality
that $n$ is divisible by $2b_n - 1$, and let $N = N(n, b_n) = n/(2b_n
-1)$ in Lemma~\ref{lemma:empiricalprocess}. Define $Y = \bar{U}' U X$.
Under the null hypothesis, $Y \sim\mathrm{N}(0, \bar{U}' \bar{U})$
and the
coordinates $Y_k$ are block-wise dependent with a bandwidth $\leq2 b_n
-1$. Split the $Y_k$'s into $2b_n-1$ different subsets $\Omega_j = \{
Y_k\dvtx k \equiv j \operatorname{mod} (2b_n -1)\}$, $1 \leq j \leq2b_n
-1$. Note that the $Y_k$'s in each subset are independent, and that
$|\Omega_j| = N$, $1 \leq j \leq2 b_n -1$.

Let $\bar{F}_n(t)$ and $\bar{F}_0(t)$ be as in the proof of Theorem
\ref
{thm:HCAdd}, and let
\[
\bar{F}_{n, j} = \frac{2b_n -1}{n} \sum_{k = 1}^n 1_{\{ Y_k^2 \geq t,
Y_k \in\Omega_j \}},\qquad 1 \leq j \leq2 b_n -1.
\]
Note that
$\bar{F}_n(t) = \frac{1}{2 b_n -1} \sum_{j = 1}^{2 b_n -1} \bar
{F}_{n,j}(t)$.
By arguments similar to that of Donoho and Jin \cite{DJ04} and basic
algebra, it follows that
\[
\mathrm{iHC}_n^* = \sup_{t} \frac{\sqrt{n} (\bar{F}_n(t) - \bar{F}_0(t))}{
\sqrt
{(2 b_n -1)\bar{F}_0(t) F_0(t)}}
\leq\sum_{j = 1}^{2b_n -1} \sup_{t} \frac{\sqrt{N} (\bar{F}_{n,j}(t)
- \bar{F}
_0(t))}{ \sqrt{\bar{F}_0(t) F_0(t)}},
\]
and so for any $x > 0$,
\[
P\{ \mathrm{iHC}_n^* \geq x \} \leq\sum_{j = 1}^{2 b_n-1} P \biggl\{
\sup
_{t} \frac{\sqrt{N} (\bar{F}_{n,j}(t) - \bar{F}_0(t))}{ \sqrt{\bar{F}_0(t)
F_0(t)}} \geq x \biggr\}.
\]
Finally, since $\bar{F}_{n,j}$'s are the empirical survival functions
of $N$
independent samples from $\chi_1^2(0)$, then
\[
\sup_{t \dvtx1/n \leq\bar{F}_0(t) \leq1/2} \frac{\sqrt{N} (\bar{F}
_{n,j}(t) -
\bar{F}_0(t))}{ \sqrt{\bar{F}_0(t) F_0(t)}} = \sup_{1/n \leq t
\leq
1/2}\mathbb{W}_N(t) \mbox{ in distribution}.
\]
Therefore,
\[
P\{ \mathrm{iHC}_n^* \geq x \} \leq(2 b_n -1) P \Bigl\{ \sup_{1/n \leq t
\leq1/2}
\mathbb{W}_N(t)\geq x \Bigr\}.
\]
Taking $x = C (\log n)^{3/2}$, the claim follows from
Lemma \ref{lemma:empiricalprocess}.
\end{pf}

\subsection{\texorpdfstring{Statement and proof of Lemma \protect
\ref
{lemma:matrixdifference}.}{Statement and proof of Lemma A.6}}

\begin{lemma} \label{lemma:matrixdifference}
The spectral norm of $\tilde{ \Sigma}_n - \bar{\Sigma}_n$ tends to zero
as $n$ tends to infinity.
\end{lemma}
\begin{pf}
To establish the lemma, note that by direct calculations and the way
$\tilde{\Sigma}$ is defined, we have
%
%
\begin{equation} \label{DefineSigma*}
\tilde{\Sigma} = \pmatrix{
\Sigma^* & \xi_{n-1} \cr
\xi_{n-1}' & 1},
\end{equation}
where
%
%
\begin{eqnarray} \label{xientries}
&&\xi_{n-1}' = \sqrt{2 n^{-\alpha}} \times\bigl(0, \ldots, n^{\alpha_0} -
k_0(n)^{\alpha},\nonumber\\[-8pt]\\[-8pt]
&&\hspace*{84.6pt} k_0(n)^{\alpha} - \bigl(k_0(n)-1\bigr)^{\alpha}, \ldots,
2^{\alpha} -1, 1\bigr),\nonumber
\end{eqnarray}
and $\Sigma^*$ is a symmetric matrix with unit diagonal entries and
with the following on the $k$th sub-diagonal:
\[
\tfrac{1}{2} \cdot\cases{
2 k^{\alpha} - (k+1)^{\alpha} - (k-1)^{\alpha}, &\quad $k \leq k_0(n)
- 1$, \cr
1 + \bigl((k-1)^{\alpha} - 2k^{\alpha}\bigr)/n^{\alpha_0} = O(n^{-\alpha
_0/\alpha
}), &\quad $k = k_0(n)$, \cr
-\bigl(1 - (k-1)^{\alpha}/n^{\alpha_0}\bigr) = O(n^{-\alpha_0/\alpha}),
&\quad
$k = k_0(n) +1$, \cr
0, &\quad $k \geq k_0(n) + 2$.}
\]
Note that $\Sigma_{n-1}(g_0)$ and $\Sigma^*$ share the $2 k_0(n) -1$
sub-diagonals that are closest to the main diagonal (including the main
diagonal). Let $H_1$ be the matrix containing all other sub-diagonals
of $\Sigma_{n-1}(g_0)$, and let $H_2$ be the matrix which contains
the $k_0(n)$th and the $(k_0(n)+1)$th diagonals (upper and lower) of
$\Sigma^*$.
It is seen that
\[
\tilde{\Sigma} - \bar{\Sigma} =
\pmatrix{
H_1 & 0 \cr
0 & 0}
+
\pmatrix{
H_2 & 0 \cr
0 & 0}
+
\pmatrix{
0 & \xi_{n-1}' \cr
\xi_{n-1}' & 0}
\equiv B_1 + B_2 + B_3.
\]
Let $\|\cdot\|_1$ and $\|\cdot\|_2$ denote the $\ell^1$ matrix norm and
the $\ell^2$ matrix norm, respectively. First, by direct calculations,
since $\alpha< 1/2$, $\|B_1 + B_2\|_1 \leq C n^{\alpha_0(\alpha-
1)/\alpha} \leq C n^{-\alpha_0}$. At the same time, by (\ref
{xientries}) and since $\alpha< 1/2$,
\[
\|B_3\|^2 \leq
\frac{C}{n^{\alpha_0}} \sum_{k = 1}^n [k^{\alpha} - (k +
1)^{\alpha
} ]^2
\leq\frac{C}{n^{\alpha_0}} \sum_{k = 1}^n k^{2 \alpha- 2} \leq
C/n^{\alpha_0}.
\]
Since the spectral norm is no greater than the $\ell^1$-matrix norm and
the $\ell^2$-matrix norm, the spectral norm of $B_1 + B_2 + B_3$ is no
greater than $C n^{-\alpha_0/2}$, and the claim follows.
\end{pf}

\subsection{\texorpdfstring{Statement and proof of Lemma
\protect\ref{lemma:LB}.}{Statement and proof of Lemma A.7}}

\begin{lemma} \label{lemma:LB}
Fix $\beta\in(\frac{1}{2}, 1)$, $r \in(0,1)$ and $\delta\in(0,1)$
such that $\bar{\gamma}_0 (1 - \delta)^{-2} r < \rho^*(\beta)$. As $n$
tends to infinity the Hellinger distance associated with model (\ref
{Model2x}) tends to zero.
\end{lemma}
\begin{pf}
To derive the lemma, let $a = \sqrt{(1 - \delta)/\bar{\gamma}_0}$, $r'
= \bar{\gamma}_0 (1 - \delta)^{-2} r $, $U_1 = a \tilde{U}$, and
$\tilde{\tilde{\mu}} = \frac{1}{a} \tilde{\mu}$. Model (\ref{Model2x})
can be
equivalently written as
%
%
\begin{equation} \label{Model4}
X = \tilde{U} \tilde{\mu} + Z = U_1 \tilde{\tilde{\mu}} + Z\qquad
\mbox
{where } Z \sim\mathrm{N}(0,I_n).
\end{equation}
Using the argument in the first paragraph of the proof of Theorem \ref
{thm:LB} it is not difficult to verify that (I) $\tilde{\tilde{\mu}}$
has $m = n^{1 - \beta}$ nonzero entries; each of which is equal to
$\sqrt{2 r' \log n}$ with $r' < \rho^*(\beta)$, and whose locations are
randomly sampled from $(1, 2, \ldots, n)$; (II) $U_1$, where $U_1(k,j)
= 0$ if $|k - j| > (\log n)^2$, is a banded lower triangular matrix and
(III) $\mathop{\overline{\lim}}_{n \rightarrow\infty}\* \max
_{\sqrt{n} \leq k \leq n -
\sqrt
{n}} (U_1' U_1) (k,k) = (1 - \delta) < 1$.

Below, write $\mu= \tilde{\tilde{\mu}}$ and $r = r'$ for short. Note
that the Hellinger distance associated with model (\ref{Model2x}) is
$E_0(\sqrt{W_n^*})$, where $E_0$ denotes the law $Z \sim\mathrm
{N}(0, I_n)$, and
\[
W_n^* = W_n^*(r, \beta; Z_1, Z_2, \ldots, Z_n) = \pmatrix{n\cr m}^{-1}
\sum_{\ell= (\ell_1, \ell_2, \ldots, \ell_m)}
e^{\mu_{\ell}' U_1' Z - \|U_1 \mu_{\ell}\|^2/2}.
\]
Introduce the set of indices
%
%
\begin{eqnarray} \label{DefineSn}
&&S_n = \Bigl\{\ell= (\ell_1, \ell_2, \ldots, \ell_m),\nonumber\\[-8pt]\\[-8pt]
&&\hspace*{28.6pt} {\min_{1
\leq j
\leq m - 1}} |\ell_{j +1} - \ell_j| \geq3 (\log n)^2, \ell_1 \geq
\sqrt
{n}, n - \ell_m \geq\sqrt{n} \Bigr\}.\nonumber
\end{eqnarray}
The following lemma is proved in Section \ref{subsub:ell}.
%
\begin{lemma} \label{lemma:ell}
Let $\ell_1 < \ell_2 < \cdots< \ell_m$ be $m$ distinct indices
randomly sampled from $(1, 2, \ldots, n)$ without replacement.
Then for any $1 \leq K \leq n$, \textup{(a)} $P\{\ell_1 \leq K\} \leq K m/n$,
\textup{(b)}
$P\{\ell_m \geq n - K\} \leq K m/n$ and \textup{(c)}\break $P\{\min_{1 \leq i \leq
m -1} \{|\ell_{i + 1} - \ell_i| \leq K \} \leq Km(m + 1)/n$. As a
result, $P\{\ell= (\ell_1, \ell_2,\break \ldots, \ell_m) \notin S_n \} =
O\{
(\log n)^2 n^{1 - 2 \beta}\} = o(1)$.
\end{lemma}

Applying Lemma \ref{lemma:ell}, we make only a negligible difference by
restricting
$\ell$ to $S_n$ and defining
%
%
\begin{equation} \label{DefineWn1}
W_n = \frac{1}{{n\choose m}} \sum_{\ell= (\ell_1, \ell_2, \ldots,
\ell
_m) \in S_n}
e^{\mu_{\ell}' U_1' Z - \|U\mu_{\ell}\|^2/2},
\end{equation}
in which case
%
%
\begin{equation} \label{decay0.32}
E(W_n^{1/2}) = E(W_n^*{}^{1/2}) + o(1).
\end{equation}
Define $Y = U_1' Z$,
%
%
\begin{equation} \label{Definesigmajs}
\sigma_j^2 = \operatorname{var}(Y_j) \equiv(U_1' U_1)(j,j),\qquad 1 \leq j
\leq n,
\end{equation}
and the event
\[
D_n = \bigl\{ Y_j / \sigma_j \leq\sqrt{2 \log n}, 1 \leq j \leq n\bigr\}.
\]
By direct calculation, $P\{D_n^c \} = o(1)$, and so by H\"older's
inequality,
\[
E\bigl(W_n^{1/2} 1_{\{D_n^c\}}\bigr) = E(W_n^{1/2}) + o(1).
\]
Combining this result and (\ref{decay0.32}) we deduce that
$E(W_n^{*{1/2}}) = E(W_n^{1/2} 1_{\{D_n\}}) + o(1)$, and comparing
this property with the desired result we see that it is is sufficient
to show that
%
%
\begin{equation} \label{decay0.4}
E\bigl(W_n^{1/2} 1_{\{D_n\} }\bigr) = 1 + o(1).
\end{equation}
The key to (\ref{decay0.4}) is the following lemma, which is proved in
Section \ref{subsub:Hellinger}.
\begin{lemma} \label{lemma:Hellinger1}
Consider the model (\ref{Model4}) where $U_1$ and $\mu$ satisfy
\textup{(I)}--\textup{(III)}. As $n \to\infty$, $E(W_n 1_{\{D_n\}}) = 1 + o(1)$,
and $E(W_n^2 1_{\{D_n\}}) = 1 + o(1)$.
\end{lemma}

Since
\[
\bigl|W_n^{1/2} 1_{\{D_n\}} - 1\bigr| \leq\frac{|W_n 1_{\{D_n\}} -1|}{1 +
W_n^{1/2} 1_{\{D_n\}}} \leq\bigl|W_n 1_{\{D_n\}} - 1\bigr|,
\]
then by H\"older's inequality,
%
%
\begin{eqnarray} \label{decay0.5}
\bigl(E\bigl|W_n^{1/2} 1_{\{D_n\}} - 1\bigr| \bigr)^2 &\leq& \bigl|W_n 1_{\{
D_n\}
} - 1 \bigr|^2 \nonumber\\[-8pt]\\[-8pt]
&=& E\bigl(W_n^2 1_{\{D_n\}} \bigr) - 2 E \bigl(W_n
1_{\{D_n\}} \bigr) + 1.\nonumber
\end{eqnarray}
Combining (\ref{decay0.5}) with Lemma \ref{lemma:Hellinger1} gives
(\ref{decay0.4}).
\end{pf}

\subsection{\texorpdfstring{Proof of Lemma \protect\ref
{lemma:ell}.}{Proof of Lemma A.8}} \label{subsub:ell}
The last claim follows once (a)--(c) are proved. Consider (a)--(b)
first. Fixing $K \geq1$, we have
\[
P\{\ell_1 = K \} = \frac{{n-K\choose m-1}}{{n\choose m}} = m \frac
{(n-m)(n - m -1) \cdots(n - m - K +2)}{n(n-1) \cdots(n - K +1)} \leq m/n,
\]
so $P\{\ell_1 \leq K\} \leq K m/n$.
Similarly,
$P\{n - \ell_m \leq K \} \leq K m/n$. This gives (a) and~(b).

Next we prove (c).
Denote the minimum inter-distance of $\ell_1, \ell_2, \ldots, \ell
_m$ by
\[
L(\ell) = L(\ell; m, n) = {\min_{1 \leq i \leq m -1}}|\ell_{i + 1} -
\ell_i|.
\]
Note that
\[
P\{L(\ell) = K\} \leq \sum_{j =1}^{m-1} P\{\ell_{j + 1} - \ell_j =
K\}
\leq\sum_{j = 1}^{m -1} \sum_{k = 1}^n P \{\ell_j = k, \ell_{j +
1} =
k + K\}.
\]
Writing
$P \{\ell_j = k, \ell_{j + 1} = k + K\} = {n\choose m}^{-1}
{k-1\choose
j -1}{n - k - K\choose m - j - 1}$,
we have
\begin{eqnarray*}
P\{L(\ell) = K\} &\leq& \frac{1}{{n\choose m}} \sum_{j =1}^{m-1} \sum_{k
= j}^n \pmatrix{k-1\cr j-1} \pmatrix{n - k - K\cr m - j -1}
\\
&=& \frac{1}{{n\choose m}} \sum_{k =1}^n \sum_{j =1}^k \pmatrix
{k-1\cr
j-1} \pmatrix{n - k - K\cr m - j -1},
\end{eqnarray*}
where the last term is no greater than
\[
\frac{1}{{n\choose m}} \sum_{k =1}^n \pmatrix{n - K -1\cr m -2}
\leq\frac{n}{{n-2\choose m}} \pmatrix{n\cr m -2} \leq m^2/n.
\]
The claim follows.

\subsection{\texorpdfstring{Proof of Lemma \protect\ref
{lemma:Hellinger1}.}{Proof of Lemma A.9}}
\label
{subsub:Hellinger}
We need the following lemma, proved in Section~\ref{proof:Truncate}.
\begin{lemma} \label{lemma:Truncate}
Consider a bivariate zero mean normal variable $(X, Y)'$ that satisfies
$\operatorname{Var}(X) = \sigma_1^2$, $\operatorname{Var}(Y) =
\sigma_2^2$ and
$\operatorname{corr}(X,Y) = \varrho$, where $c_0 \leq\sigma_1,
\sigma_2
\leq
1$ for some constant $c_0 \in(0,1)$. Then there is a constant $C > 0$
such that, for sufficiently large $n$,
\begin{eqnarray*}
E \bigl[\exp( A_n X - \sigma_1^2 A_n^2/2 ) \cdot1_{\{Y >
\sigma_2 T_n\}} \bigr]& \leq& C \cdot n^{- (1 - \varrho\sqrt{r})^2}\\
&\leq& C n^{-(1 - \sqrt{r})^2},
\\
E \biggl[\exp\biggl( A_n(X + Y) - \frac{\sigma_1^2 + \sigma_2^2}{2}A_n^2
\biggr) \cdot1_{\{X \leq\sigma_1 T_n, Y \leq\sigma_2 T_n\}} \biggr]
&\leq& C n^{-d(r)},
\end{eqnarray*}
where $d(r) = \min\{2r, 1 - 2 (1 - \sqrt{r})^2\}$.
\end{lemma}

Now we proceed with the derivation of Lemma \ref{lemma:Hellinger1}.
Consider the first claim. Note that for any $\ell= (\ell_1, \ell_2,
\ldots, \ell_m) \in S_n$, the minimum inter-distance of $\ell_i$ is no
less than $3 (\log n)^2$, and so
\[
\|U_1 \mu_{\ell}\|^2 = A_n^2 \sum_{i = 1}^m (U_1' U_1)(\ell_i, \ell_i)
= A_n^2 \sum_{i = 1}^m \sigma_{\ell_i}^2.
\]
In view of the definition of $Y_j$ and $\sigma_j$ [see (\ref
{Definesigmajs})], we can rewrite $W_n$ as
%
%
\begin{equation} \label{hel1}
W_n = \frac{1}{{n\choose m}} \sum_{\ell= (\ell_1, \ell_2, \ldots,
\ell
_m) \in S_n} \exp\Biggl( A_n \sum_{i = 1}^m Y_{\ell_i} - \frac
{A_n^2}{2} \sum_{i =1}^m \sigma_{\ell_i}^2 \Biggr).
\end{equation}
Note that
%
%
\begin{equation} \label{hel2}
1_{\{D_n^c\}} \leq\sum_{j =1}^{n} 1_{ \{Y_j / \sigma_j > T_n \} }.
\end{equation}
Combining (\ref{hel1}) and (\ref{hel2}) gives
%
%
\begin{eqnarray} \label{Toshow1.5}
&&E\bigl(W_n \cdot1_{\{D_n^c\}}\bigr) \nonumber\\
&&\qquad\leq\frac{1}{{n\choose m}} \sum_{\ell=
(\ell_1, \ldots, \ell_m) \in S_n} \sum_{k =1}^{n} E \Biggl[ \exp
\Biggl(
A_n \sum_{j =1}^m Y_{\ell_j} - \frac{A_n^2}{2} \sum_{j = 1}^m
\sigma
_{\ell_j}^2 \Biggr) \\
&&\qquad\quad\hspace*{186.7pt}{}\times 1_{\{Y_k/\sigma_k > T_n\}} \Biggr].\nonumber
\end{eqnarray}
We shall say that two indices $j$ and $k$ are near each other if $|j -
k| \leq(\log n)^2$. In this notation, for each $1 \leq k \leq n$, when
$k$ is near one $\ell_j$, say $\ell_{j_0}$, $Y_k$ must be independent
of all other $Y_{\ell_j}$ with $j \neq j_0$. It follows that
\begin{eqnarray*}
&&E \Biggl[\exp\Biggl( A_n \sum_{j =1}^m Y_{\ell_j} - \frac{A_n^2}{2}
\sum
_{j = 1}^m \sigma_{\ell_j}^2 \Biggr) \cdot1_{\{Y_k/\sigma_k > T_n\}}
\Biggr]\\
&&\qquad = E \bigl[ \exp( A_n Y_{\ell_{j_0}} - \sigma_{j_0}^2
A_n^2/2 ) \cdot1_{\{Y_k / \sigma_k > T_n\}} \bigr].
\end{eqnarray*}
By Lemma \ref{lemma:Truncate}, the right-hand side is no greater than
$C n^{ - (1 - \sqrt{r})^2}$. Therefore,
%
%
\begin{equation} \label{Toshow2}
\qquad E \Biggl[\exp\Biggl( A_n \sum_{j =1}^m Y_{\ell_j} - \frac{A_n^2}{2}
\sum
_{j = 1}^m \sigma_{\ell_j}^2 \Biggr) \cdot1_{\{Y_k/\sigma_k > T_n\}}
\Biggr] \leq C n^{-(1 - \sqrt{r})^2}.
\end{equation}
Moreover, for each fixed $\ell= (\ell_1, \ldots, \ell_m) \in S_n$,
there are at most $2 m (\log n)^2$ different indices $k$ that can be
near some of the $\ell_j$'s; and when they are, they can be near only
one such $\ell_j$.
Combining these results gives
%
%
\begin{eqnarray} \label{hel3}
E\bigl[W_n \cdot1_{\{D_n^c\}}\bigr] &\leq& \frac{1}{{n\choose m}} \sum_{\ell=
(\ell_1, \ldots, \ell_m) \in S_n} C (\log n)^2 m n^{-(1 -
\sqrt{r})^2}\nonumber\\[-8pt]\\[-8pt]
&\leq& C (\log n)^2 n^{(1 - \beta) - (1 - \sqrt{r})^2}.\nonumber
\end{eqnarray}
By the definition of $\rho^*(\beta)$ and the assumption of the lemma,
$r < \rho^*(\beta) \leq(1 - \sqrt{1 - \beta})^2$, and so
the first claim follows directly from (\ref{hel3}).

We now consider the second claim. Fix $0 \leq N \leq m$, and
let $\tilde{S}_N(\ell)$ denote the set of all $k = (k_1, k_2, \ldots,
k_m) \in S_n$ such that there are exactly $N$ $k_j$'s that are near to
one $\ell_i$. (Clearly, any $k_j$ can be near to at most one
$\ell_i$.) The two sets of indices $(\ell_1, \ell_2, \ldots, \ell_m)$
and $(k_1, k_2, \ldots, k_m)$ form exactly $N$ pairs where each
contains one candidate from the first set and one candidate from the second.
These pairs are not near to each other and not near to any remaining
indices outside the pairs.
Using (\ref{hel1}), we write
%
%
\begin{eqnarray} \label{HC}
&&E\bigl[W_n^2 \cdot1_{\{D_n\}}\bigr] \nonumber\\
&&\qquad= \pmatrix{n\cr m}^{-2}\sum_{\ell= (\ell_1, \ell_2, \ldots, \ell_m) \in
S_n}\nonumber\\[-0pt]\\[-8pt]
&&\qquad\quad\hspace*{33.2pt}{}\times \sum_{N = 0}^m
\sum_{k = (k_1, k_2, \ldots, k_m) \in\tilde{S}_N(\ell)}
E \Biggl[ \exp\Biggl(A_n \sum_{i = 1}^m (Y_{\ell_i} +
Y_{k_i})\nonumber\\
&&\hspace*{258pt}\hspace*{-40.8pt}{} - \frac{A_n^2}{2} \sum_{i =1}^m (\sigma_{\ell_i}^2 +
\sigma
_{k_i}^2) \Biggr) \cdot1_{\{D_n\}} \Biggr].\nonumber
\end{eqnarray}
For any fixed $\ell$ and $k \in\tilde{S}_N(\ell)$, by symmetry, and
without loss of generality, we suppose the $N$ pairs are $(\ell_1,
k_1)$, $(\ell_2, k_2), \ldots, (\ell_N, k_N)$. By independence of the
pairs with other indices, and also by independence among the pairs,
%
%
\begin{eqnarray}\label{Toshow9.1}
&& E \Biggl[ \exp\Biggl( A_n \sum_{j =1}^m (Y_{\ell_j} + Y_{k_j}) -
\frac
{A_n^2}{2} \sum_{j = 1}^m (\sigma_{\ell_j}^2 + \sigma_{k_j}^2)
\Biggr)
\cdot1_{\{D_n\}} \Biggr] \nonumber\\
&&\qquad\leq E \Biggl[ \exp\Biggl( A_n \sum_{j =1}^m (Y_{\ell_j} + Y_{k_j}) -
\frac{A_n^2}{2} \sum_{j = 1}^m (\sigma_{\ell_j}^2 + \sigma_{k_j}^2 )
\Biggr) \nonumber\\
&&\qquad\quad\hspace*{46.8pt}{}\times1_{\{Y_{\ell_j}/\sigma_{\ell_j} \leq T_n, Y_{k_j} /
\sigma
_{k_j} \leq T_n,\ \mathrm{for}\ \mathrm{all}\ 1 \leq j \leq N \}} \Biggr]
\nonumber\\
&&\qquad\leq E \Biggl[ \exp\Biggl( A_n \Biggl\{\sum_{j =1}^{N} (Y_{\ell_j} +
Y_{k_j}) - \frac{A_n^2}{2}\sum_{j = 1}^{N} (\sigma_{\ell_j}^2
+\sigma
_{k_j}^2 ) \Biggr\} \Biggr) \\
&&\qquad\quad\hspace*{58pt}{}\times1_{\{Y_{\ell_j}/\sigma_{\ell_j}
\leq
T_n, Y_{k_j}/\sigma_{k_j} \leq T_n,\ \mathrm{for}\ \mathrm{all}\ 1 \leq j \leq
N\}} \Biggr] \nonumber\\
&&\qquad= \prod_{j = 1}^{N} \biggl( E \biggl[ \exp\biggl\{ A_n (Y_{\ell_j} +
Y_{k_j}) - \frac{A_n^2}{2} (\sigma_{\ell_j}^2 + \sigma_{k_j}^2)
\biggr\}
\nonumber\\
&&\qquad\quad\hspace*{90.1pt}{}\times1_{\{Y_{\ell_j}/\sigma_{\ell_j} \leq T_n, Y_{k_j}/\sigma_{k_j}
\leq T_n\}} \biggr] \biggr).\nonumber
\end{eqnarray}
Here, in the first inequality, we have used the fact that
\[
1_{\{D_n\}} \leq1_{\{Y_{\ell_j}/\sigma_{\ell_j} \leq T_n,
Y_{k_j}/\sigma_{k_j} \leq T_n,\ \mathrm{for}\ \mathrm{all}\ 1 \leq j \leq N\}};
\]
in the second inequality, we have utilized the independence and the
fact that
\[
E[ \exp( A_n Y_j - \sigma_j^2 A_n^2/2 )] = 1\qquad
\mbox{for all } j = 1, \ldots, n,
\]
and in the third equality, we have used again the independence.
Moreover, in view of the definition of $U_1$, and Lemma \ref
{lemma:Udecay}, there is a constant $c_0 \in(0,1)$ such that
$\sigma_j \in[c_0,1]$.
Using Lemma \ref{lemma:Truncate}, for sufficiently large $n$ and each
$1 \leq j \leq N$,
%
%
\begin{eqnarray} \label{Toshow9}
&&E \biggl[ \exp\biggl( A_n (Y_{\ell_j} + Y_{k_j}) - \frac{A_n^2}{2}
(\sigma
_{\ell_j}^2 + \sigma_{k_j}^2) \biggr) \nonumber\\[-8pt]\\[-8pt]
&&\hspace*{70.6pt}{}\times 1_{\{Y_{\ell_j}/\sigma
_{\ell
_j} \leq T_n, Y_{k_j}/\sigma_{k_j} \leq T_n\}} \biggr] \leq
C n^{d(r)}\nonumber
\end{eqnarray}
with $d(r)$ being as in Lemma \ref{lemma:Truncate}.
Combining (\ref{Toshow9.1}) and (\ref{Toshow9}) gives
%
%
\begin{equation} \label{Toshow10.5}
E\bigl[W_n^2 \cdot1_{\{D_n\}}\bigr] \leq\pmatrix{n\cr m}^{-2} \sum_{\ell=
(\ell_1, \ldots, \ell_m)} \sum_{N = 0}^m \bigl(C n^{d(r)}\bigr)^N
|\tilde
{S}_{N}(\ell) |,
\end{equation}
where $|\tilde{S}_N(\ell)|$ denotes the cardinality of $\tilde
{S}_N(\ell
)$. By elementary
combinatorics,
%
%
\begin{eqnarray} \label{H6}
|\tilde{S}_N(\ell)| &\leq& \pmatrix{m\cr N} (2 \log^2 n )^N \pmatrix
{n -
N\cr m - N} \nonumber\\[-8pt]\\[-8pt]
&\leq& (2 \log^2 n )^N \pmatrix{m\cr N} \pmatrix{n\cr m - N}.\nonumber
\end{eqnarray}
Direct calculations show that
%
%
\begin{equation} \label{H7}
\frac{{m\choose N} {n\choose m - N}}{{n\choose m}} = \frac{1}{N!}
\biggl(\frac{m!}{(m-N)!} \biggr)^2 \frac{(n-m)!}{(n- m + N)!} \lesssim
\frac
{1}{N!} \biggl(\frac{m^2}{n} \biggr)^N.
\end{equation}
Substituting (\ref{H6}) and (\ref{H7}) into (\ref{Toshow10.5}) and
recalling that $m = n^{1 - \beta}$, we deduce that
%
%
\begin{eqnarray} \label{Toshow11}\qquad
&&E\bigl[W_n^2 \cdot1_{\{D_n\}}\bigr]\nonumber\\[-8pt]\\[-8pt]
&&\qquad \leq\pmatrix{n\cr m}^{-1} \sum_{\ell=
(\ell_1, \ell_2, \ldots, \ell_m) \in S_n} \sum_{N = 0}^m \frac{1}{N!}
\biggl(\frac{m^2}{n} \biggr)^N \bigl( (C \log^2 n) n^{ d(r)}
\bigr)^N,\nonumber
\end{eqnarray}
where the last term does not exceed $\sum_{N = 0}^{\infty
}(N!)^{-1}
[C(\log^2 n) n^{1 + d(r) - 2 \beta} ]^N$.
By the assumption of the lemma,
\[
r < \rho^*(\beta) =
\cases{
\beta- 1/2, &\quad $1/2 < \beta\leq3/4$, \cr
\bigl(1 - \sqrt{1 - \beta}\bigr)^2, &\quad $3/4 \leq\beta< 1$;}
\]
thus it can be seen that $1 + d(r) - 2 \beta< 0$ for all fixed $\beta$
and $r \in(0, \rho^*(\beta))$. Combining this with (\ref{Toshow11})
gives the second claim.

\subsection{\texorpdfstring{Proof of Lemma \protect\ref
{lemma:Truncate}.}{Proof of Lemma A.10}}
\label
{proof:Truncate}
Denote the density, cdf and survival function of $\mathrm{N}(0,1)$ by
$\phi$,
$\Phi$ and $\bar{\Phi}$.
For the first claim, define $W = X/\sigma_1$ and $V = Y/\sigma_2$ if
$\rho\geq0$ and $V = - Y/\sigma_2$ otherwise. The proofs for two
cases $\rho\geq0$ and $\rho< 0$ are similar, so we only show the
first one. In this case, it suffices to show that
\[
E \bigl[\exp( \sigma_1 A_n W - \sigma_1^2 A_n^2/2 )
\cdot
1_{\{V > T_n\}} \bigr] \leq C \cdot n^{- (1 - \varrho\sqrt{r})^2}.
\]
Write $W = (W - \rho V) + \rho V$, and note that $(1 - \rho)^2 + \rho^2
\leq1$. It is seen that
%
%
\begin{eqnarray} \label{Truncate1}
\sigma_1 A_n W - \sigma_1^2 A_n^2/2 &\leq& [\sigma_1 A_n (W - \rho V) -
\sigma_1^2 (1 - \rho)^2 A_n^2/2]\nonumber\\[-8pt]\\[-8pt]
&&{} + [\sigma_1 A_n \rho V - \sigma_1^2
\rho^2 A_n^2/2].\nonumber
\end{eqnarray}
Since $W$ and $V$ have unit variance and correlation $\rho$, then $(W -
\rho V)$ is independent of $V$ and is distributed as $\mathrm{N}(0, (1
- \rho
)^2)$. Therefore, $E[\exp( \sigma_1 A_n (W - \rho V) - \sigma_1^2 (1 -
\rho)^2 A_n^2/2)] = 1$. Combining this with (\ref{Truncate1}) gives
\begin{eqnarray*}
&&E \bigl[\exp( \sigma_1 A_n W - \sigma_1^2 A_n^2/2 )
\cdot
1_{\{V > T_n\}} \bigr] \\
&&\qquad= E \bigl[\exp( \sigma_1 \rho A_n V -
\sigma_1^2 \rho^2 A_n^2/2 ) \cdot1_{\{V > T_n\}} \bigr].
\end{eqnarray*}
Now, by direct calculation,
\[
E \bigl[\exp( A_n V - A_n^2/2 ) \cdot1_{\{V > T_n\}}
\bigr]
= \int_{T_n}^{\infty} \phi(x - \sigma_1 \rho A_n) \,dx = \bar{\Phi
}(T_n -
\sigma_1 \rho A_n).
\]
Since $\bar{\Phi}(x) \leq C \phi(x)$ for all $x > 0$,
\[
\bar{\Phi}(T_n - \sigma_1 \rho A_n) \leq C \phi(T_n - \sigma_1
\rho A_n)
= C
n^{- (1 - \rho\sqrt{r})^2}.
\]
Combining these results gives the claim.

We now establish the second claim. By H\"older's inequality, it
suffices to show that
\[
E\bigl[\exp(2 A_n X - \sigma_1^2 A_n^2 ) \cdot1_{\{ X \leq\sigma_1 T_n\}
}\bigr] \leq C n^{-d(r)}.
\]
Recalling that $W = X / \sigma_1$, we have
\[
E\bigl[\exp(2 A_n X - \sigma_1 A_n^2 ) \cdot1_{\{ X \leq\sigma_1 T_n\}}\bigr]
= E\bigl[\exp(2 \sigma_1 A_n W - \sigma_1^2 A_n^2 ) \cdot1_{\{ W \leq
T_n\}}\bigr].
\]
By direct calculation,
\begin{eqnarray*}
E\bigl[\exp(2 \sigma_1 A_n W - \sigma_1^2 A_n^2 ) \cdot1_{\{ W \leq T_n\}
}\bigr] &=& e^{\sigma_1^2 A_n^2} \int_{-\infty}^{T_n} \phi(x - 2 \sigma_1 A_n)
\,dx \\
&=& e^{\sigma_1^2 A_n^2} \Phi(T_n - 2 \sigma_1 A_n).
\end{eqnarray*}
Since $\Phi(x) \leq C \phi(x)$ for all $x < 0$ and $\Phi(x) \leq1$ for
all $x \geq0$,
\begin{eqnarray*}
&&e^{\sigma_1^2 A_n^2} \Phi(T_n - 2 \sigma_1 A_n) \\
&&\qquad\leq\cases{
C e^{\sigma_1^2 A_n^2} = C n^{2 \sigma_1^2 r}, &\quad $\sigma_1^2 r
\leq
1/4$, \cr
e^{\sigma_1^2 A_n^2} \phi(T_n -2 \sigma_1 A_n) = C n^{1 - 2 (1 -
\sigma
_1 \sqrt{r})^2}, &\quad $\sigma_1^2 r > 1/4$.}
\end{eqnarray*}
In view of the definition of $d(r)$, $e^{\sigma_1^2 A_n^2} \Phi(T_n - 2
\sigma_1 A_n) \leq C n^{d (\sigma_1^2 r)}$. Since that $\sigma_1
\leq
1$ and that $d(r)$ is a monotonely increasing function, we have
$d(\sigma_1^2 r) \leq d(r)$. Combining these results gives the claim.

\mbox{}

\subsection{\texorpdfstring{Statement and proof of Lemma
\protect\ref{lemma:HCproof}.}{Statement and proof of Lemma A.11}}

\begin{lemma} \label{lemma:HCproof}
Under the conditions of Theorem \ref{thm:HC}, the right-hand side of
(\ref{HCmodelxx}) converges to zero algebraically fast as $n$
diverges to infinity.
\end{lemma}
\begin{pf}
The key observation needed to establish the lemma is that there is a
sequence of positive numbers $\delta_n$ that tends to zero as $n$
diverges to infinity such that $\nu_k \geq(1 - \delta_n) A_n$ for all
$k \in\{\ell_1, \ell_2, \ldots, \ell_m\}$, so it is natural to compare
model (\ref{HCmodelx1}) with the following model:
%
%
\begin{equation} \label{HCmodelx2}
Y^* = \nu^* + Z,\qquad Z \sim\mathrm{N}(0,I_n),
\end{equation}
where $\nu^*$ has $m$ nonzero entries of equal strength $(1 - \delta_n)
A_n$ whose locations are randomly drawn from $\{1, 2, \ldots, n\}$
without replacement.

For short, write $t = t_n^*$ and
\[
h_n(t) = \frac{\sqrt{n}(\bar{F}_n(t) - \bar{F}_0(t))}{\sqrt{(2 b_n
-1) \bar{F}_0(t)
(1 - \bar{F}_0(t))}}.
\]
Let $\bar{F}_n^*(t)$ be the empirical survival function of $\{
(Y_k^*)^2\}
_{k = 1}^n$, and let $\bar{F}(t) = E[\bar{F}_n(t)]$ and $\bar
{F}^*(t) = E[\bar{F}
_n^*(t)]$. Recall that the family of noncentral $\chi^2$-distributions
has monotone likelihood ratio. Then $\bar{F}(t) \geq\bar{F}^*(t)
\geq\bar{F}
_0(t)$. Now, first, since the $Y_k$'s are block-wise dependent with a
block size $\leq2b_n - 1$, it follows by direct calculations that
\[
\operatorname{Var}(h_n(t)) \leq C \bar{F}(t)/\bar{F}_0(t).
\]
Second, by $\bar{F}(t) \geq\bar{F}_n^*(t)$,
%
%
\begin{eqnarray} \label{HCmodelxx2}
E[h_n(t)] &=& \frac{\sqrt{n}( \bar{F}(t) - \bar{F}_0(t))}{\sqrt{(2
b_n -1)
\bar{F}
_0(t) (1 - \bar{F}_0(t))}} \nonumber\\[-8pt]\\[-8pt]
&\geq&\frac{\sqrt{n}( \bar{F}^*(t) - \bar{F}
_0(t))}{\sqrt
{(2 b_n -1) \bar{F}_0(t) (1 - \bar{F}_0(t))}},\nonumber
\end{eqnarray}
where the right-hand side diverges to infinity algebraically fast by an
argument similar to that in \cite{DJ04}. Combining Chebyshev's
inequality, the identity $b_n = \log n$ and calculations of the mean
and variance of $h_n(t)$, we deduce that
%
%
\begin{equation} \label{HCmodelxx1}
P\{ h_n(t) \leq(\log n)^2\} \leq C (\log n) \frac{\bar{F}(t)}{n
(\bar{F}
(t) - \bar{F}_0(t) )^2}.
\end{equation}

It remains to show that the last term in (\ref{HCmodelxx1}) is
algebraically small.
We discuss separately the cases $\bar{F}(t)/\bar{F}_0(t) \geq2$ and
$\bar{F}
(t)/\bar{F}
_0(t) < 2$. For the first case,
\[
\frac{\bar{F}(t)}{n(\bar{F}(t) - \bar{F}_0(t))^2} \leq\frac{C}{n
\bar{F}(t)}
\leq
\frac{C}{n \bar{F}_0(t)},
\]
which is algebraically small since $t = \sqrt{2 q \log n}$ and $0 < q <
1$. For the second case,
%
%
\begin{equation} \label{HCmodelxx3}
\frac{\bar{F}(t)}{n(\bar{F}(t) - \bar{F}_0(t))^2} \leq\frac{C
\bar{F}_0(t)}{n
(\bar{F}(t)
- \bar{F}_0(t))^2} \leq\frac{C \bar{F}_0(t)}{n (\bar{F}^*(t) -
\bar{F}_0(t))^2},
\end{equation}
which is seen to be algebraically small by comparing it to the
right-hand side of (\ref{HCmodelxx2}).
\end{pf}

\subsection{\texorpdfstring{Statement and proof of Lemma
\protect\ref{lemma:PD}.}{Statement and proof of Lemma A.12}}

\begin{lemma} \label{lemma:PD}
Let $\Sigma_n$ be as in (\ref{hjmatrix}). For sufficiently large $n$,
necessary and sufficient conditions for $\Sigma_n$ to be positive
definite are, respectively, $0\leq\alpha\leq2$ and $0 < \alpha_0
\leq
\alpha\leq1$.
\end{lemma}
\begin{pf}
We begin by establishing the first claim. Suppose such an
autoregressive structure exists for $\alpha\geq\alpha_0 > 0$. Let
\[
Y_k = \sqrt{a_n} \cdot(X_{k + 1} - X_k)/d,\qquad a_n = n^{\alpha
_0}/2, k = 1, 2, \ldots, n-1.
\]
Clearly, $\operatorname{var}(Y_k) = 1$. At the same time, direct calculation
shows that
the correlation between $Y_1$ and $Y_{j+1}$ equals to
$[(j + 1)^{\alpha} + (j - 1)^{\alpha} - 2 j^{\alpha}]/2$ for all $1
\leq j \leq n-2$, which is no larger than $1$. Taking $j = 2$ yields
$(3^{\alpha} + 1 - 2 \cdot2^{\alpha})/2 \leq1$,
and hence $\alpha\leq2$.

Consider the second claim. For any $k \geq1$, define the partial sum
$S_{k}(t) = 1 + 2 \sum_{j = 1}^{k} (1 - \frac{j^{\alpha}}{n^{\alpha
_0}})^+ \cos(kt)$. By a well-known result in trigonometry \cite
{Zygmund}, to establish the positive-definiteness of $\Sigma_n$, it
suffices to show that
%
%
\begin{eqnarray} \label{Eq1.1}
S_{k_0 + 1}(t) &\geq& 0 \qquad\mbox{for all $t \in[-\pi, \pi]$}\quad
\mbox{and}\nonumber\\[-8pt]\\[-8pt]
S_{k_0+1}(t) &>& 0\qquad\mbox{except for a set of measure zero}.\nonumber
\end{eqnarray}
Here, $k_0 = k_0(n; \alpha, \alpha_0)$ is the largest integer $k$
such that
$k^{\alpha} \leq n^{\alpha_0}$.

We now derive (\ref{Eq1.1}).
Using a result from \cite{Zygmund}, page 183, if we let
$a_0 = 2$, and $a_j = 2 (1 - \frac{j^{\alpha}}{n^{\alpha_0}})^+$, $1
\leq j \leq n-1$, then
$S_{k_0 + 1}(t) = \sum_{j = 0}^{k_0 -1} (j + 1) \Delta^2 a_{j} K_j(t) +
(k_0 + 1) K_{k_0}(t) \Delta a_{k_0} + D_n(t) a_{k_0 + 1}$.
Here,
$\Delta a_j = a_{j} - a_{j + 1}$, $\Delta^2 a_j = a_{j} + a_{j + 2} - 2
a_{j+1}$, and
$D_j(t)$ and $K_j(t)$ are the Dirichlet's kernel and the Fej\'er's
kernel, respectively,
%
%
\begin{eqnarray} \label{DK}
D_j(t) &=& \frac{\sin((j + {1/2})t) }{2
\sin({t/2})},\nonumber\\[-8pt]\\[-8pt]
K_j(t) &=& \frac{2}{j+1} \biggl(\frac{\sin(({j+1})/{2}t)}{2 \sin
({t/2})} \biggr)^2,\qquad j = 0, 1, \ldots.\nonumber
\end{eqnarray}
In view of the definition of $k_0$,
$a_{k_0 + 1} = (1 - \frac{(k_0 + 1)^{\alpha}}{n^{\alpha_0}})^+ = 0$. Also,
by the monotonicity of $\{a_j\}$,
$\Delta a_{k_0} = a_{k_0} - a_{k_0 + 1} \geq0$.
Therefore,
$S_{k_0 + 1}(t) \geq\sum_{j = 0}^{k_0 -1} (j + 1) \Delta^2 a_{j} K_j(t)$.

We claim that the sequence $\{a_0, a_1, \ldots, a_{n-1}\}$ is convex.
In detail, since $\alpha\leq1$, the sequence $\{j^{\alpha}\}$ is
concave. As a result, the sequence $\{(1 - \frac{j^{\alpha
}}{n^{\alpha
_0}})\}$ is convex, and so is the sequence $\{(1 - j^{\alpha
}/n^{\alpha
_0})^+ \}$. In view of the definition of $a_j$, the claim follows
directly. The convexity of the $a_j$'s implies that
$\Delta^2 a_j \geq0$, $0 \leq j \leq n-2$. Therefore, $S_{k_0 +
1}(t) \geq0$.
This proves the first part of (\ref{Eq1.1}).

We now prove the second part of (\ref{Eq1.1}), and discuss separately
the two cases $\alpha< 1$ and $\alpha= 1$. In the first case,
$\Delta a_0 = n^{-\alpha_0}(2 - 2^{\alpha}) > 0$ and $K_0(t) = \frac
{1}{2}$. As a result, $S_{k_0 + 1}(t) \geq(2 - 2^{\alpha})/(2
n^{\alpha_0}) > 0$, and the claim follows. In the second case,
$\Delta a_j = n^{-\alpha_0}(2 j - j - (j+2)) = 0$,
and
$\Delta a_{k_0 - 1} = [1 - n^{-\alpha_0}(k_0-1)] - 2 (1 - n^{-\alpha
_0}k_0) = n^{-\alpha_0} (k_0 + 1) -1 > 0$. Therefore,
$S_{k_0 + 1}(t) \geq(k_0 + 1) [(k_0 + 1)n^{-\alpha_0} - 1] K_{k_0}(t)$.
Clearly, $S_{k_0 + 1}(t)$ can only assume $0$ when $\frac{1}{2}(k_0 +
1) t$ is a multiple of $\pi$. Since the set of such $t$ has measure
zero, the claim follows directly.
\end{pf}

\subsection{\texorpdfstring{Statement and proof of Lemma
\protect\ref{lemma:f0}.}{Statement and proof of Lemma A.13}}

\begin{lemma} \label{lemma:f0}
For $0 < \alpha< 1$, we have $\mathrm{essinf}_{- \pi\leq\theta\leq
\pi} \{f_{\alpha}(\theta)\} > 0$.
\end{lemma}
\begin{pf}
To derive the lemma, let $a_0 = 2$, and $a_k = 2 k^{\alpha} -
(k+1)^{\alpha} - (k-1)^{\alpha}$, $1 \leq k \leq n-1$.
Clearly, $a_k > 0$ for all $k$, so $f_{\alpha}(0; \alpha) > 0$.
Furthermore, when $\theta\neq0$,
by \cite{Zygmund}, equation 1.7, page 183,
%
%
\begin{equation} \label{zygmund1}
f_{\alpha}(\theta) = \sum_{\nu= 0}^{\infty} (\nu+1) [a_{\nu+2} +
a_{\nu} - 2 a_{\nu+ 1}] a_{\nu} K_{\nu}(\theta),
\end{equation}
where $K_{\nu}(\theta)$ is the Fej\'er's kernel as in (\ref{DK}). By
the positiveness of the Fej\'er's kernel, all remains to show is that
$a_{k+1} + a_{k-1} -2a_k > 0$, for all $k \geq2$.

Define
$h(x) = (1 + 2x)^{\alpha} + (1 - 2x)^{\alpha} - 4(1 + x)^{\alpha} - 4(1
- x)^{\alpha} + 6$, $0 \leq x \leq1/2$.
By direct calculations, for all $k \geq2$,
%
%
\begin{eqnarray} \label{Eq1.1a}\hspace*{32pt}
&&a_{k+1} + a_{k-1} -2a_k \nonumber\\
&&\qquad= - k^{\alpha} \biggl[ \biggl(1 + \frac
{2}{k}
\biggr)^{\alpha}
+ \biggl(1 - \frac{2}{k} \biggr)^{\alpha} - 4 \biggl(1 + \frac{1}{k}
\biggr)^{\alpha}
- 4 \biggl(1 - \frac{1}{k} \biggr)^{\alpha} + 6 \biggr] \\
&&\qquad= -k^{\alpha}h(1/k).\nonumber
\end{eqnarray}
Also, by basic calculus,
\[
h''(x) = 4\alpha(\alpha-1) [(1 + 2x)^{\alpha-2} + (1 - 2x)^{\alpha-2}
- (1+x)^{\alpha-2} - (1 - x)^{\alpha-2}].
\]
Since $0 < \alpha< 1$, $x^{\alpha-2}$ is a convex function. It
follows that
$h''(x) < 0$ for all $x \in(0, 1/2)$, and $h(x)$ is a strictly
concave function. At the same time, note that $h(0) = h'(0) = 0$, so
$h(x) < 0$ for $x \in(0, 1/2]$. Combining this with (\ref{Eq1.1a})
gives the claim.
\end{pf}
\end{appendix}

\section*{\texorpdfstring{Acknowledgments.}{Acknowledgments}}
Jiashun Jin would like to thank Christopher Genovese and Larry
Wasserman for extensive discussion, and Aad van der Vaart for help on
the proof of (\ref{monotone}). He would also like to thank Peter Bickel, Emmanuel
Cand\'es, David Donoho, Karlheinz Gr\"ochenig, Michael Leinert,
Joel Tropp and Zepu Zhang for encouragement and pointers.


%
\printaddresses

\end{document}